\documentclass[preprint,10t]{elsarticle}

\usepackage{amssymb}
\usepackage{stmaryrd}
\usepackage{tipa}
\usepackage{amsfonts,amsmath,latexsym}
\usepackage{mathrsfs}
\usepackage{amsbsy}
\usepackage{indentfirst}
\usepackage{amsmath}
\usepackage{graphicx,psfrag,epsf}
\usepackage{enumerate}
\usepackage{natbib}
\usepackage{url} 
\usepackage{amsfonts}
\usepackage{mathrsfs}
\usepackage{algorithm}
\usepackage{algorithmic}
\usepackage{subfigure}
\usepackage{epsfig}
\usepackage{booktabs}
\usepackage{multirow}

\numberwithin{equation}{section}

\newtheorem{theorem}{\textbf{Theorem}}[section]
\newtheorem{lemma}{\textbf{Lemma}}

\newtheorem{ass}{\textbf{Assumption}}

\newtheorem{remark}[theorem]{Remark}

\newcommand{\beqnar}{\begin{eqnarray*}}
\newcommand{\eeqnar}{\end{eqnarray*}}
\newcommand{\ba}{\begin{array}}
\newcommand{\ea}{\end{array}}

\newenvironment{proof}[1]{\begin{trivlist}\item {\it
\bf Proof.}\quad} {\qed\end{trivlist}}

\allowdisplaybreaks[4]

\journal{}

\begin{document}

\begin{frontmatter}



\title{ Bias Correction Estimation for Continuous-Time Asset Return Model with Jumps }
\author{Yuping Song\footnote{Corresponding author, email: songyuping@shnu.edu.cn }}

\author{Ying Chen}

\author{Zhouwei Wang}

\address{ School of Finance and Business, Shanghai Normal
         University, Shanghai, 200234, P.R.China}

\begin{abstract}
In this paper, local linear estimators are adapted for the unknown
infinitesimal coefficients associated with continuous-time asset
return model with jumps, which can correct the bias automatically
due to their simple bias representation. The integrated diffusion
models with jumps, especially infinite activity jumps are mainly
investigated. In addition, under mild conditions, the weak
consistency and asymptotic normality is provided through the
conditional Lindeberg theorem. Furthermore, our method presents
advantages in bias correction through simulation whether jumps
belong to the finite activity case or infinite activity case.
Finally, the estimators are illustrated empirically through the
returns for stock index under five-minute high sampling frequency
for real application.
\end{abstract}


\end{frontmatter}

\section{Introduction}

Continuous-time models are widely used in economics and finance,
such as interest rate etc, especially the continuous-time diffusion
processes with jumps. Jump-diffusion process $X_{t}$ is represented
by the following stochastic differential equation:
\begin{equation}
\label{1.1} dX_{t} = \mu(X_{t-})dt + \sigma(X_{t-})dW_{t} +
\int_{\mathscr{E}}c(X_{t-} , z) r(\omega, dt, dz), \end{equation}
which can accommodate the impact of sudden and large shocks to
financial markets. Johannes \cite{joh} provided the statistical and
economic role of jumps in continuous-time interest rate models.
However, in empirical finance the current observation usually
behaves as the cumulation of all past perturbation such as stock
prices by means of asset returns in Nicolau \cite{n2} et al.
Furthermore, in the research field involved with model (\ref{1.1}),
although the scholars focused on the price of the asset, most of
them didn't consider the returns of the asset. As mentioned in
Campbell, Lo and MacKinlay \cite{clm}, return series of an asset are
a complete and scale-free summary of the investment opportunity for
most investors, and are easier to handle than price series due to
their more attractive statistical properties. The capital asset
return model is the core of capital market theory which describes
the relationship between the returns and risks of individual
securities or portfolio.

For characterizing this integrated economic phenomenon, moreover,
the return series, we considered the promising continuous-time
integrated diffusion process with jumps (\ref{om}), which is
motivated by unit root processes under the discrete framework of
Park and Phillips \cite{pp} and continuous integrated diffusion
process in Nicolau \cite{n}. It satisfies the following second-order
stochastic differential equation:
\begin{equation}
\label{om} \left\{
\begin{array}{ll}
dY_{t} = X_{t}dt,\\
dX_{t} = \mu(X_{t-})dt + \sigma(X_{t-})dW_{t}+
\int_{\mathscr{E}}c(X_{t-} , z) r(\omega, dt, dz),
\end{array}
\right.\end{equation} where $\mathscr{E} =
\mathbb{R}\setminus\{0\},$ $W=\{W_t\}_{t\ge 0}$ is a standard
Brownian motion, $r(\omega, dt, dz) = (p - q)(dt, dz), ~p(dt, dz)$
is a time-homogeneous Poisson random measure on $\mathbb{R}_{+}
\times \mathbb{R}$, which is independent of $W_{t}$, and $q(dt, dz)$
is its intensity measure, that is, $E[p(dt, dz)] = q(dt, dz) =
f(z)dzdt$, $f(z)$ is a L\'{e}vy density. For empirical financial
data, $X_{t}$ in the model (\ref{om}) represents the continuously
compounded return of underlying assets, $Y_{t}$ denotes the asset
price by means of the cumulation of the returns plus initial asset
value. Furthermore, the model (\ref{om}) can accommodate
nonstationarity and transform nonstationarity into stationarity by
differencing, which can not be performed through univariate
diffusion model due to the nondifferentiability of a Brownian
motion.

For model (\ref{om}), the estimators for unknown coefficients have
been considered based on low frequency or high frequency
observations under various settings. For model (\ref{om}) without
jumps, Gloter \cite{g1}\cite{g2} and Ditlevsen and S{\o}rensen
\cite{ds} built the parametric and semiparametric estimation, while
Nicolau \cite{n} and Comte, Genon-Catalot and Rozenholc \cite{cgr}
analyzed nonparametric estimators for the unknown quantities.
Moreover, Wang and Lin \cite{wl}, Wang, Zhang and Tang \cite{wzt},
Hanif \cite{hm2} and Wang and Tang \cite{wt} improved the
nonparametric estimators for them. For model (\ref{om}) with finite
activity jumps, Song and Lin \cite{sl}, Song, Lin and Wang
\cite{slw}, Chen and Zhang \cite{cz} and Funke and Schmisser
\cite{fs} theoretically investigated nonparametric estimation for
the drift or volatility coefficients. Song \cite{sy} empirically
considered the application for the estimation proposed in high
frequency financial data.

In this paper, we adapt local linear estimators for the unknown
coefficients of integrated diffusion models with jumps, especially
infinite activity jumps. In the context of nonparametric estimator
with finite-dimensional auxiliary variables, local polynomial
smoothing become an effective smoothing method, which doesn't assume
the functional form for the unknown coefficients. Moreover, local
linear estimators have excellent properties such as full asymptotic
minimax efficiency achievement and boundary bias correction
automatically, one can refer to Fan and Gijbels \cite{fg2} for
better review.

Our contribution have three folds. Firstly, in terms of the model,
the previous work was mainly focused on the continuous case as
Nicolau \cite{n} or finite activity case as Song \cite{sy}. We will
consider a more practical integrated diffusion models with infinite
activity jumps for the asset return. The existence of infinite
activity jumps for the high frequency financial data has been
testified based on A\"it-Sahalia and Jacod \cite{aj} in empirical
analysis part.

Secondly, in the theoretical side, compared with Nadaraya-Watson
estimators, the conditional Lindeberg theorem might be no longer
applicable for the local linear estimators due to their destroyed
adaptive and predictable structure of conditionally on the
$\sigma-$field generated by $X_{t}.$ We effectively tackle the key
technical problems by means of the Slutsky's theorem and establish
central limit theorems for the volatility functions in the
second-order diffusion model with infinite activity jumps. More
technical proof details can be sketched in Lemma \ref{l4} and
Theorem \ref{mr}.

Thirdly, Chen and Zhang \cite{cz} gave the large sample properties
of local linear estimators for second-order diffusion with finite
activity jump, but they didn't consider the finite-sampling
performance of them. Considering what has been talked above, in the
practical side we consider two types of jump (finite activity jumps
and infinite activity jumps) aimed at verifying the better
finite-sampling performance of local linear estimators under various
settings. Moreover, the estimators are illustrated empirically
through the return of stock index in Shenzhen Stock Exchange under
five-minute high sampling frequency data between Jan 2015 and Dec
2015. In summary, the integrated diffusion model with jump,
especially infinite activity jumps may be an alternative model to
describe the dynamic variation for the returns of financial assets.

The paper is organized as follows. The local linear estimators and
their large sample properties are collected in Section 2. The finite
sample performance of underlying estimators through Monte Carlo
simulation study is presented in Section 3. The estimators are
illustrated empirically in Section 4. Some technical lemmas for the
main theorems are given in Appendix part.

\section{ Local Linear estimators and Large sample properties }

For model (\ref{om}), we usually get observations
$\{Y_{i\Delta_n};i=1,2,\cdot \cdot \cdot\}$ rather than
$\{X_{i\Delta_n};i = 1, 2, \cdot \cdot \cdot\}.$ However, the value
of $X_{t_{i}}$ cannot be obtained from $Y_{t_{i}} = Y_{0} +
\int_{0}^{t_{i}} X_{s} ds$ in a fixed sample intervals.
Additionally, nonparametric estimations of the unknown qualities in
model (\ref{om}) cannot in principle be constructed on the
observations $\{Y_{i\Delta_{n}} ; i = 1, 2, \cdot \cdot \cdot\}$ due
to the unknown conditional distribution of $Y$. As Nicolau \cite{n}
showed, with observations $\{Y_{i\Delta_n};i=1,2,\cdot \cdot
\cdot\}$ and given that
$$Y_{{i\Delta_{n}}} - Y_{{(i-1)\Delta_{n}}} = \int_{(i-1)\Delta_{n}}^{i\Delta_{n}} X_{u} du, $$
we can obtain an approximation value of $X_{i\Delta_n}$ by
\begin{equation}
\label{yi}
\widetilde{X}_{i\Delta_n}=\frac{Y_{i\Delta_n}-Y_{(i-1)\Delta_n}}{\Delta_n}.
\end{equation}

Due to the Markov properties of model (\ref{om}), we can build the
following infinitesimal conditional expectations
\begin{align}
& \label{ice1} E\left[\frac{\widetilde{X}_{(i+1)\Delta_{n}} -
\widetilde{X}_{i\Delta_{n}}}{\Delta_{n}}|\mathscr{F}_{(i-1)\Delta_{n}}\right]
 = \mu (X_{(i-1)\Delta_{n}}) + O_{p}(\Delta_{n}), \\
& \label{ice2} E\left[\frac{(\widetilde{X}_{(i+1)\Delta_{n}} -
\widetilde{X}_{i\Delta_{n}})^{2}}{\Delta_{n}}|\mathscr{F}_{(i-1)\Delta_{n}}\right]
= \frac{2}{3}\sigma^{2} (X_{(i-1)\Delta_{n}}) +
\frac{2}{3}\int_{\mathbb{R}}c^{2}(X_{(i-1)\Delta_{n}} , z)f(z)dz +
O_{p}(\Delta_{n}).
\end{align}
where $\mathscr{F}_t=\sigma\{X_s,s \leq t\}$. One can refer to
Appendix A in Song, Lin and Wang \cite{slw} for detailed
calculations.

For the given $\{\widetilde{X}_{i\Delta_{n}} ; i= 1, 2, \cdots\}$,
the local linear estimators for $\mu(x)$ and $M(x)$ based on
infinitesimal conditional expectations (\ref{ice1}) and (\ref{ice2})
are defined as the solutions to the following weighted least squares
problems: find $a_1$, $b_1$, $a_2$, $b_2$ to minimize
\begin{equation}
\label{ll1} \sum_{i=1}^{n}\Big(\frac{\widetilde{X}_{(i+1)\Delta_{n}}
- \widetilde{X}_{i\Delta_{n}}}{\Delta_{n}}
 - a_{1} - b_{1}\big({\widetilde{X}_{i\Delta_{n}} - x}\big) \Big)^{2}
  K\Big(\frac{\widetilde{X}_{(i-1)\Delta_{n}} - x}{h_{n}}\Big), \end{equation}
\begin{equation}
\label{ll2}
\sum_{i=1}^{n}\Big(\frac{\frac{3}{2}(\widetilde{X}_{(i+1)\Delta_{n}}
- \widetilde{X}_{i\Delta_{n}})^{2}}{\Delta_{n}}
 - a_{2} - b_{2}\big({\widetilde{X}_{i\Delta_{n}} - x}\big) \Big)^{2} K\Big(\frac{\widetilde{X}_{(i-1)\Delta_{n}} - x}{h_{n}}\Big),
\end{equation}
where $K(\cdot)$ is the kernel function and $h_{n}$ is a sequence of
positive numbers, satisfies $h_n\to 0$ as $n\to \infty.$

The solutions for $a_{1}$ and $a_{2}$ to (\ref{ll1}) and (\ref{ll2})
as follows are respectively the local linear estimators of $\mu(x)$
and $M(x) = \sigma^2(x)+\int_{\mathscr{E}}c(x,z)f(z)dz,$
\begin{equation}
\label{ll3} \hat{\mu}_n(x) = \frac{\sum_{i=1}^n \omega_{i-1}\Big(
\frac{\widetilde{X}_{i+1}-\widetilde{X}_{i}}
{\Delta_n}\Big)}{\sum_{i=1}^n \omega_{i-1}},
\end{equation}
\begin{equation}
\label{ll4} \hat{M}_{n}(x) = \frac{\sum_{i=1}^n \omega_{i-1}
\frac{3}{2}
\frac{(\widetilde{X}_{i+1}-\widetilde{X}_{i})^2}{\Delta_n}}{\sum_{i=1}^n
\omega_{i-1}}
\end{equation}
where
\begin{equation*}\omega_{i-1} = K\left(\frac{\widetilde{X}_{i-1} - x}{h_{n}}\right)\Bigg(\sum_{j=1}^n K\left(\frac{\widetilde{X}_{j-1} - x}{h_{n}}\right)
(\widetilde{X}_{j} - x)^2 - (\widetilde{X}_{i} - x) \sum_{j=1}^{n}
K\left(\frac{\widetilde{X}_{j-1} - x}{h_{n}}\right) (\widetilde{X}_j
- x)\Bigg). \end{equation*}

The assumptions of this paper are listed below, which confirm the
large sample properties of the constructed estimators based on
(\ref{ll3}) and (\ref{ll4}).

\begin{ass}\label{a1}
i) (Local Lipschitz continuity)~ For each $n \in \mathbb{N},$ there
exist a constant $L_{n}$ and a function $\zeta_{n}$ : $\mathscr{E}
\rightarrow \mathbb{R}_{+}$ with $\int_{\mathscr{E}}
\zeta_{n}^{2}(z) \nu(dz) < \infty$ such that, for any $|x| \leq n,
|y| \leq n, z \in {\mathscr{E}}$,
$$|\mu(x) - \mu(y)| + |\sigma(x) - \sigma(y)| \leq L_{n}|x - y|,~~~ |c(x , z) - c(y , z)| \leq \zeta_{n}(z)|x - y|.$$

(ii (Linear growthness)~ For each $n \in \mathbb{N}$, there exist
$\zeta_{n}$ as above and C, such that for all $x \in \mathbb{R}, z
\in {\mathscr{E}}$,
$$|\mu(x)| + |\sigma(x)| \leq C (1 + |x|) , ~ |c(x , z)| \leq \zeta_{n}(z)(1 + |x|).$$\end{ass}

\begin{remark}
\label{ar1} This assumption guarantees the existence and uniqueness
of a solution to stochastic differential equation $X_{t}$ in
(\ref{om}), see Jacod and Shiryaev \cite{js}. For instance, Long, Ma
and Shimizu \cite{lms} and Long and Qian \cite{lq} imposed similar
conditions on the coefficients of the underlying stochastic
differential equation.
\end{remark}

\begin{ass} \label{a2}
The process $X=\{X_{t}\}_{\ge 0}$ is ergodic and stationary with a
finite invariant measure $\phi(x)$. For a given point $x_{0}$, the
stationary probability measure $p(x)$ of the process $X$ is positive
at $x_{0},$ that is $p(x_{0})
> 0.$ Furthermore, the process $X$ is $\rho-$mixing with
$\sum_{i\geq1}\rho(i\Delta_{n}) = O(\frac{1}{\Delta_{n}^{\alpha}}),
~n \rightarrow \infty,$ where $\alpha < 1/2.$
\end{ass}

\begin{remark}
\label{ar2} The Assumption \ref{a2} implies that the process $X_{t}$
has a unique weak solution. The finite invariant measure implies
that the process $X_{t}$ is positive Harris recurrent with the
stationary probability measure $p(x) =
\frac{\phi(x)}{\phi(\mathscr{D})}, ~\forall x \in \mathscr{D}.$ The
hypothesis that $X_{t}$ is a stationary process is obviously a
plausible assumption because for major integrated time series data,
a simple differentiation generally assures stationarity. The same
condition yielding information on the rate of decay of $\rho-$mixing
coefficients for $X_{t}$ was mentioned the Assumption 3 in
Gugushvili and Spereij \cite{gs}. For instance, the $\rho-$mixing
process $X_{t}$ with exponentially decreasing mixing coefficients
satisfies the condition, see Hansen and Scheinkman \cite{hs}, Chen,
Hansen and Carrasco \cite{chc}.
\end{remark}

\begin{ass} \label{a3}
The kernel $K$($\cdot$) : $\mathbb{R} \rightarrow \mathbb{R}^{+}$ is
a positive and continuously differentiable function satisfying:
$$\int K(u) du = 1,~K_{i}^{j} := \int K^{i}(u) u^{j} du < \infty.$$
Moreover, For $2\le i\le n$,
$$\lim_{h\rightarrow0}E\left[~\frac{1}{h}~|K^{'}(\xi_{n,i})|^{\alpha}
\Big(\frac{\widetilde{X}_{(i-1)\Delta_{n}} - x}{h}\Big)^{m}\right] <
\infty$$ where $\alpha= 1, 2$ or $ 4$, $m = 0, 1$ or $2$ and
$\xi_{n,i}$ = $\theta\left(\frac{X_{(i-1)\Delta_{n}} - x}{h}\right)
+ (1 - \theta)\left(\frac{\widetilde{X}_{(i-1)\Delta_{n}} -
x}{h}\right)$, $0 \leq \theta \leq 1.$
\end{ass}

\begin{remark}
\label{ar3} In fact, any density function can be considered as a
kernel, moreover even unnecessary positive functions can be used.
For simplification, we only consider positive and symmetrical
kernels used widely. It is well known both empirically and
theoretically that the choice of kernel functions is not very
important to the kernel estimator, see Gasser and M\"{u}ller
\cite{gm}. As Nicolau \cite{n} pointed out this assumption is
generally satisfied under very weak conditions. For instance, with a
Gaussian kernel and a Cauchy stationary density (which has heavy
tails) we still have
$\lim_{h\rightarrow0}E\big[(\frac{1}{h})|K^{'}(\frac{X}{h})|^{4}\big]
< \infty$. Notice that the expectation with respect to the
distribution $\xi_{n,i}$ depends on the stationary densities of
$X_{n,i}$ and $\widetilde{X}_{n,i}$ because $\xi_{n,i}$ is a convex
linear combination of $X_{n,i}$ and $\widetilde{X}_{n,i}.$
\end{remark}

\begin{ass} \label{a4}
For every $p \geq 1$, $\sup_{t\geq0} E[|X_{t}|^{p}] < \infty,$ and
$\int_{\mathscr{E}}|z|^{p}\nu(dz) < \infty.$
\end{ass}

\begin{remark}
\label{ar4} This assumption guarantees that Lemma \ref{l1} can be
used properly throughout the article. If $X$ is a L\'{e}vy process
with bounded jumps (i.e., $\sup_{t}|\Delta X_{t}| \leq C < \infty$
almost surely, where C is a nonrandom constant), then
$E\{|X_{t}^{n}|\} < \infty ~\forall n$, that is, $X_{t}$ has bounded
moments of all orders, see Protter \cite{pr}. This condition is
widely used in the estimation of an ergodic diffusion or
jump-diffusion from discrete observations, see Florens-Zmirou
\cite{fz}, Kessler \cite{ke}, Shimizu and Yoshida \cite{syo}.
\end{remark}

\begin{ass} \label{a5}
$\Delta_{n} \rightarrow 0, ~
(\frac{n\Delta_{n}}{h_{n}})(\Delta_{n}\log(\frac{1}{\Delta_{n}}))^{\frac{1}{2}}
\rightarrow 0, ~h_{n}n\Delta_{n}^{1+\alpha} \rightarrow \infty,
h_{n}^{5}n\Delta_{n} \rightarrow 0~as~ n \rightarrow \infty.$
\end{ass}

\begin{remark}
\label{ar5} The relationship between $h_{n}$ and $\Delta_{n}$ is
similar as the stationary case in Bandi and Nguyen \cite{bn} and
(b1), (b2) of A8 in Nicolau \cite{n}.
\end{remark}

We have the following asymptotic results for the local linear
estimators such as (\ref{ll3}) and (\ref{ll4}) based on the
assumptions above.

\begin{theorem}
\label{mr} Under Assumptions \ref{a1}-\ref{a5}, as $n \rightarrow
\infty,$ we have

(i)~~$\hat{\mu}_{n}(x) \stackrel{P} \longrightarrow
\mu(x),~\hat{M}_{n}(x) \stackrel{P} \longrightarrow M(x).$

(ii)~Furthermore, if $h_{n} n \Delta_{n} \longrightarrow 0$ and
$h_{n} = O((n \Delta_{n})^{-1/5}),$ then
$$\sqrt{h_{n} n \Delta_n}\left(\hat{\mu}_{n} (x) - \mu(x) - \frac{1}{2}h_{n}^{2}\mu^{''}(x)\frac{(K_{1}^{2})^{2} -
K_{1}^{3}K_{1}^{1}}{K_{1}^{2} - (K_{1}^{1})^{2}}\right)
\Longrightarrow N\left(0, V\frac{M(x)}{p(x)}\right),$$ and
$$\sqrt{h_{n} n \Delta_n}\left(\hat{M}_{n} (x) - M(x) - \frac{1}{2}h_{n}^{2}M^{''}(x)\frac{(K_{1}^{2})^{2} -
K_{1}^{3}K_{1}^{1}}{K_{1}^{2} - (K_{1}^{1})^{2}}\right)
\Longrightarrow
N\left(0,V\frac{\int_{\mathscr{E}}c^4(x,y)f(y)dy}{p(x)}\right),$$
where
$$V=\frac{(K_1^2)^2 K_2^0+ (K_1^1)^2 K_2^2 -2(K_1^1)(K_1^2)K_2^1}{[K_1^2-(K_1^1)^2]^2}.$$
\end{theorem}

\begin{remark}
\label{mrr1} For the finite activity jumps of model (\ref{om}),
$$\int_{\mathscr{E}}c(X_{t-} , z) r(\omega, dt, dz) :=
\int_{\mathscr{E}} c(X_{t-}, z) N(dt, dz) -
\lambda(X_{t-})\int_{\mathscr{E}} c(X_{t-}, z)\Pi(dz)dt,$$ where
$N(dt, dz)$ is a Poisson counting measure, $c(\cdot, y)$ reflects
the conditional impact of a jump and $\Pi(dz)$ is the probability
distribution function of a jump. For model (\ref{om}), we can
observe that $f(z) = \lambda(x) \Pi(z),$ so $M(x) = \sigma^{2}(x) +
\lambda(x)\int_{\mathscr{E}} c^{2}(x, z)\Pi(dz)$ and
$\int_{\mathscr{E}}c^4(x,z)f(z)dz = \lambda(x)\int_{\mathscr{E}}
c^{4}(x, z)\Pi(dz).$

For the infinite activity jumps of model (\ref{om}), we will focus
on the following diffusion process with jumps for $X_{t}$ such that
\begin{equation*}
d X_{t} = \mu(X_{t-}) d t + \sigma(X_{t-}) d W_{t} + \xi(X_{t-}) d
J_{t},
\end{equation*}
where $J_{t}$ is a pure jump L\'{e}vy process such as an infinite
activity jump process of the representation
$$d J_{t} = \int_{\mathbb{R}} y (\mu(dt, dy) - \nu(dy)dt) := \int_{\mathbb{R}}y \bar{\mu}(dt,
dy)$$ with $\mu(dt, dy)$ a Poisson random measure compensated by its
intensity measure $\nu(dy)dt.$ Hence $M(x) = \sigma^{2}(x) +
\int_{\mathscr{R}} c^{2}(x, z)\nu(dz)$ and
$\int_{\mathscr{R}}c^4(x,z)f(z)dz = \int_{\mathscr{R}} c^{4}(x,
z)\nu(dz)$ for model (\ref{om}).

In contrary to the integrated diffusion model without jumps (Nicolau
\cite{n}), the rate of convergence of the second infinitesimal
moment estimator is same as the first infinitesimal moment
estimator. Apparently, this is due to the presence of discontinuous
breaks that have an equal impact on all the functional estimates. As
Johannes \cite{joh} pointed out, for the conditional variance of
interest rate changes, not only diffusion play a certain role, but
also jumps account for more than half at lower interest level rates,
almost two-thirds at higher interest level rates, which dominate the
conditional volatility of interest rate changes. Thus, it is
extremely important to estimate the conditional variance as
$\sigma^{2}(x)$ + $\int_{\mathscr{E}}c^{2}(x , z)f(z)dz$ which
reflects the fluctuation of the return of the underlying asset.
\end{remark}

\begin{remark}
\label{mrr2} In Song \cite{sy}, he showed that under the Assumptions
in this paper, $h_{n}n\Delta_{n} \rightarrow 0$ and $h_{n} =
O((n\Delta_{n})^{-1/5}),$ the following result holds with symmetric
kernels
$$\sqrt{h_{n}n\Delta_{n}}\left(\hat{\mu}^{NW}_{n}(x) - \mu(x) - h_{n}^{2}
K_{1}^{2} \left(\frac{1}{2}\mu^{''}(x) +
\mu^{'}(x)\frac{\phi^{'}(x)}{\phi(x)}\right)\right)
\stackrel{d}{\rightarrow} N\Big( 0,
K_{2}^{0}\frac{M(x)}{p(x)}\Big).$$ When $K$ is symmetric, for local
linear estimators we obtain that $$\sqrt{h_{n} n
\Delta_n}\left(\hat{\mu}_{n} (x) - \mu(x) -
\frac{1}{2}h_{n}^{2}\mu^{''}(x)K_{1}^{2}\right) \Longrightarrow
N\left(0, K_{2}^{0}\frac{M(x)}{p(x)}\right).$$ Comparing the bias of
local linear estimator with that of Nadaraya-Watson estimator in
Song \cite{sy} above, we can observe that the bias of
$\hat{\mu}_{n}^{NW} - \mu(x)$ is
$\mu^{'}(x)\frac{\phi^{'}(x)}{\phi(x)}K_{1}^{2}$ more than that of
$\hat{\mu}_{n}(x) - \mu(x)$. When $K$ is asymmetric, sketching the
proof procedure in Song \cite{sy}, we can prove that the bias for
$\hat{\mu}^{NW} - \mu(x)$ should subtract $hK_{1}^{1},$ which
confirms the bias of local linear estimator is smaller than that of
Nadaraya-Watson case. Hence, compared with the Nadaraya-Watson
estimator, local linear estimator possesses simple bias
representation and can correct the bias automatically whether
$K(\cdot)$ is symmetric or not.
\end{remark}

\begin{remark}
\label{mrr3} Similarly as Theorem 4 in Fan and Gijbels \cite{fg1},
we now investigate the behavior of the estimator (\ref{ll3}) at
left-boundary points (the right-boundary points are the same). Put
$x_{n} = ch_{n}$, with $c
> 0$. Assuming $nh_{n} \rightarrow \infty$, the conditional MSE of
the estimator (\ref{ll3}) at the boundary point $x_{n}$ is
$O(h_{n}^{4} + \frac{1}{nh_{n}})$. More proof details see Fan and
Gijbels \cite{fg1} (A little tedious, so we omit here). Indeed, its
rate of convergence for local linear estimator is not influenced by
the position of the point under consideration. Hence the local
linear smoother does not require modifications at the boundary. So
it turns out that the local linear smoother has an additional
advantage over other kernel-type estimator (see Gasser and
M\"{u}ller \cite{gm}). To some extend, the simulation in the paper
confirms this result.
\end{remark}

\begin{remark}
\label{mrr4} It is very important to consider the choice of the
bandwidth in nonparametric estimation. Here we will select the
optimal bandwidth $h_{n}$ based on the mean squared error (MSE) and
the asymptotic theory in Theorem \ref{mr}. Take $\mu(x)$ for
example, the optimal smoothing parameter $h_{n}$ for local linear
estimator of $\mu(x)$ is given that
\begin{equation*}
h_{n,opt,fi} = \left(\frac{1}{n\Delta_{n}} \cdot
\frac{4VM(x)(K_{1}^{2} -
(K_{1}^{1})^{2})^{2}}{p(x)\mu^{2''}(x)((K_{1}^{2})^{2} -
K_{1}^{3}K_{1}^{1})^{2}} \right)^{\frac{1}{5}}
 = O_{p}\left(\frac{1}{n\Delta_{n}}\right)^{\frac{1}{5}},\end{equation*}
which differs from the continuous case with
$h_{n,opt}=O_{p}({n}^{-1/5}).$ The bandwidth $h_{n}$ constructed
above relies on the consistent estimators for these unknown
quantities and they are difficult to obtain and may give rise to
bias. Here we mention two rules of thumb on selecting the bandwidth.
Practically, for simplicity one can use the empirical bandwidth
selector $h_{n} = 1.06*\hat{S}*(n\Delta_{n})^{-1/5},$ where
$\hat{S}$ denotes the standard deviation of the data. Or, one can
apply the cross-validation method to assess the performance of an
estimator via estimating its prediction error. The main idea is to
minimize the following expression: $CV(h) =
n^{-1}\sum_{i=1}^{n}\{\frac{\widetilde{X}_{(i+1)\Delta_{n}} -
\widetilde{X}_{i\Delta_{n}}}{\Delta_{n}} -
\hat{a}_{h,-i}(\widetilde{X}_{i\Delta_{n}})\}^{2}$, where
$\hat{a}_{h,-i}(\widetilde{X}_{i\Delta_{n}})$ is the local linear
estimator (\ref{ll3}) and bandwidth $h$, but without using the $i$th
observation.
\end{remark}

\begin{remark}
\label{mrr5} If the smoothing parameter $h_{n} =
O((n\Delta_{n})^{-1/5}),$ the normal confidence interval for
$\mu(x)$ using local linear estimators at the significance level
$100(1-\alpha)\%$ are constructed as follows,
\begin{align*}
I_{\mu,\alpha} = & \Bigg[\hat{\mu}_{n}(x) - h_{n}^{2} \cdot
\frac{1}{2} \hat{\mu}_{n}^{''}(x)\frac{(K_{1}^{2})^{2} -
K_{1}^{3}K_{1}^{1}}{K_{1}^{2} - (K_{1}^{1})^{2}} - z_{1-\alpha/2}
\cdot \frac{1}{\sqrt{n\Delta_{n} h_{n}}} \cdot
\sqrt{V \frac{\hat{M}_{n}(x)}{\hat{p}_{n}(x)}},\\
& \hat{\mu}_{n}(x) - h_{n}^{2} \cdot \frac{1}{2}
\hat{\mu}_{n}^{''}(x)\frac{(K_{1}^{2})^{2} -
K_{1}^{3}K_{1}^{1}}{K_{1}^{2} - (K_{1}^{1})^{2}} + z_{1-\alpha/2}
\cdot \frac{1}{\sqrt{n\Delta_{n} h_{n}}} \cdot \sqrt{V
\frac{\hat{M}_{n}(x)}{\hat{p}_{n}(x)}} \Bigg],
\end{align*}
where $z_{1-\alpha/2}$ is the inverse CDF for the standard normal
distribution evaluated at $1 - \alpha/2.$ To facilitate statistical
inference for $\mu(x)$ based on Theorem \ref{mr}, we need to conduct
consistent estimators for the unknown quantities in the normal
approximation. $\hat{\mu}_{n}(x),$ $\hat{M}_{n}(x)$ denote the local
linear estimators of $\mu(x), M(x)$ in (\ref{ll3}) and (\ref{ll4}),
respectively. As Fan and Gijbels \cite{fg2} showed, the derivative
$\hat{\mu}_{n}^{''}(x)$ in $I_{\mu,\alpha}$ can be estimated by
taking the second derivative of the local linear estimators of
$\mu(x)$ in (\ref{ll3}). The consistent estimator for $p(x)$ is
$$\hat{p}_{n}(x) =
\frac{1}{nh_{n}} \sum_{i=1}^{n} K
\left(\frac{\widetilde{X}_{(i-1)\Delta_{n}} - x}{h_{n}}\right).$$
\end{remark}

\section{ Monte Carlo Simulation Study }
In this section, a simple Monte Carlo simulation experiment is
constructed aimed at the finite-sampling performance between local
linear estimators constructed as (\ref{ll3}) and (\ref{ll4}) denoted
as LL and Nadaraya-Watson estimators constructed in Song (2017)
denoted as NW. Throughout this section, various lengths of
observation time interval $T~(= 50, 100, 500)$ and sample sizes
$n~(= 500, 1000, 2500)$ with $\Delta_{n} = \frac{T}{n}$ will be
considered. We use classical Gaussian kernel $K(x) =
\frac{1}{\sqrt{2\pi}}e^{-\frac{x^{2}}{2}}$ and the common bandwidth
$h = 1.06*\hat{S}*(n\Delta_{n})^{-\frac{1}{5}} =
1.06*\hat{S}*T^{-\frac{1}{5}}$, where $\hat{S}$ denotes the standard
deviation of the data. Here we will consider two types of jump aimed
at verifying the better finite-sampling performance of local linear
estimators under various settings. Take the estimator
$\hat{\mu}_{n}(x)$ of $\mu(x)$ for example to show the
small-sampling performance. Similar results can be found for
$\hat{M}_{n}(x)$ of $M(x).$

{\textbf{Example 1 (Finite Activity case).}} Our experiment is based
on the following second-order diffusion model with finite activity
jump:
\begin{equation}
\label{mc1} \left\{
\begin{array}{ll}
dY_{t} = X_{t-}dt,\\
dX_{t} = -10X_{t-}dt + \sqrt{0.1+0.1X_{t-}^{2}}dW_{t} + dJ_{t},
\end{array}
\right. \end{equation} where the coefficients of continuous part are
equal to the ones used in Nicolau \cite{n} and $J_{t}$ is a compound
Poisson jump process with arrival intensity $\lambda = 2$ and jump
size $Z_{n} \sim \mathscr{N}(0,0.036^{2})$ corresponding to Bandi
and Nguyen \cite{bn}.

The finite-sample performance of the local linear estimator and NW
estimator for $\mu(x)$ with finite activity jump is demonstrated in
Figure 1. We can observe that local linear estimator performs a
little better than the NW estimator, especially at the boundary
points. What is shown in Table 1 is the biases of local linear
estimator and NW estimator at various quantile points of $X_{t}.$ In
addition, 95\% Monte Carlo confidence intervals of various
estimators for $\mu(x)$ are depicted in Figure 2. It shows that the
true value of $\mu(x)$ can not fall in the 95\% Monte Carlo
confidence intervals of NW estimators at the sparse design boundary
points. These findings confirms the fact that local linear estimator
can correct the boundary bias automatically due to its simple bias
representation as shown in Remark \ref{mrr2}. Figure 3 gives the QQ
Plots for local linear estimator of the drift function $\mu(x)$ with
finite activity jump, which reveals the normality of local linear
estimators in finite sample and confirms the results in Theorems
\ref{mr}.

\begin{figure}[!htb]
\centering
\includegraphics[width=1 \textwidth]{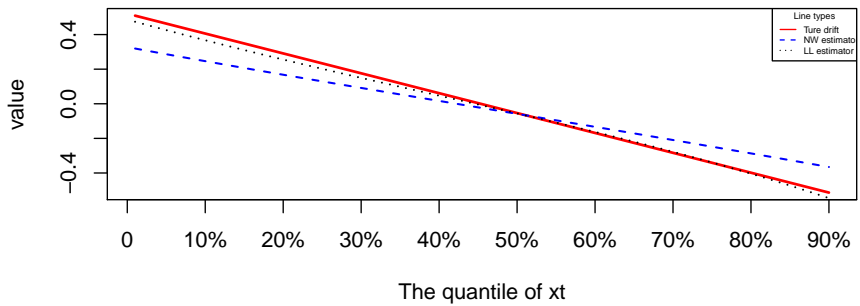}
\caption{ Various Nonparametric Estimators for $\mu(x) = - 10*x$
with $T = 10,~n = 1000,~\lambda = 2$ and jump size $Z_{n} \sim
\mathscr{N}(0,0.036^{2})$}
\end{figure}

\begin{table}[!htb]
  \centering
  \caption{The Biases of Nadaraya-Watson and Local Linear Estimators for $\mu(x) = -
  10*x$ at various quantile points of sample $X_t$ with $T = 10,~n = 1000,~\lambda = 2$ and jump size $Z_{n} \sim \mathscr{N}(0,0.036^{2})$ }
    \begin{tabular}{cc|c|c|c|c|c|c|c|c}
    \toprule
    \multirow{2}[4]{*}{\textbf{Bias}} & \multicolumn{9}{c}{\textbf{Various Quantile Points of Sample Xt}} \\
\cmidrule{2-10}          & \multicolumn{1}{c}{\textbf{10\%}} & \multicolumn{1}{c}{\textbf{20\%}} & \multicolumn{1}{c}{\textbf{30\%}} & \multicolumn{1}{c}{\textbf{40\%}} & \multicolumn{1}{c}{\textbf{50\%}} & \multicolumn{1}{c}{\textbf{60\%}} & \multicolumn{1}{c}{\textbf{70\%}} & \multicolumn{1}{c}{\textbf{80\%}} & \textbf{90\%} \\
    \midrule
    \textbf{NW} & \multicolumn{1}{c}{-0.1898 } & \multicolumn{1}{c}{-0.1302 } & \multicolumn{1}{c}{-0.0817 } & \multicolumn{1}{c}{-0.0468 } & \multicolumn{1}{c}{-0.0193 } & \multicolumn{1}{c}{0.0178 } & \multicolumn{1}{c}{0.0556 } & \multicolumn{1}{c}{0.0965 } & 0.1465  \\
    \textbf{LL} & \multicolumn{1}{c}{-0.0668 } & \multicolumn{1}{c}{-0.0388 } & \multicolumn{1}{c}{-0.0232 } & \multicolumn{1}{c}{-0.0161 } & \multicolumn{1}{c}{-0.0131 } & \multicolumn{1}{c}{-0.0128 } & \multicolumn{1}{c}{-0.0178 } & \multicolumn{1}{c}{-0.0305 } & -0.0600  \\
    \bottomrule
    \end{tabular}%
  \label{tab:addlabel}%
\end{table}%

\begin{figure}[!htb]
\centering
\includegraphics[width=1 \textwidth]{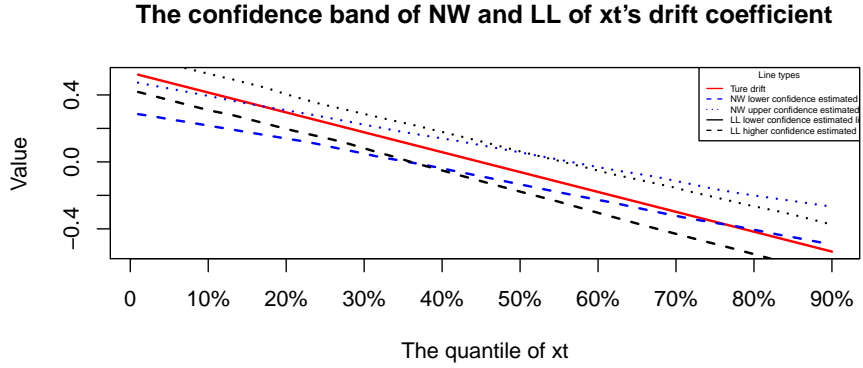}
\caption{ 95\% Monte Carlo Confidence Interval for $\mu(x) = - 10*x$
with $T = 10,~n = 1000,~\lambda = 2$ and jump size $Z_{n} \sim
\mathscr{N}(0,0.036^{2})$ }
\end{figure}

\begin{figure}[!htb]
\centering
\includegraphics[width=1 \textwidth]{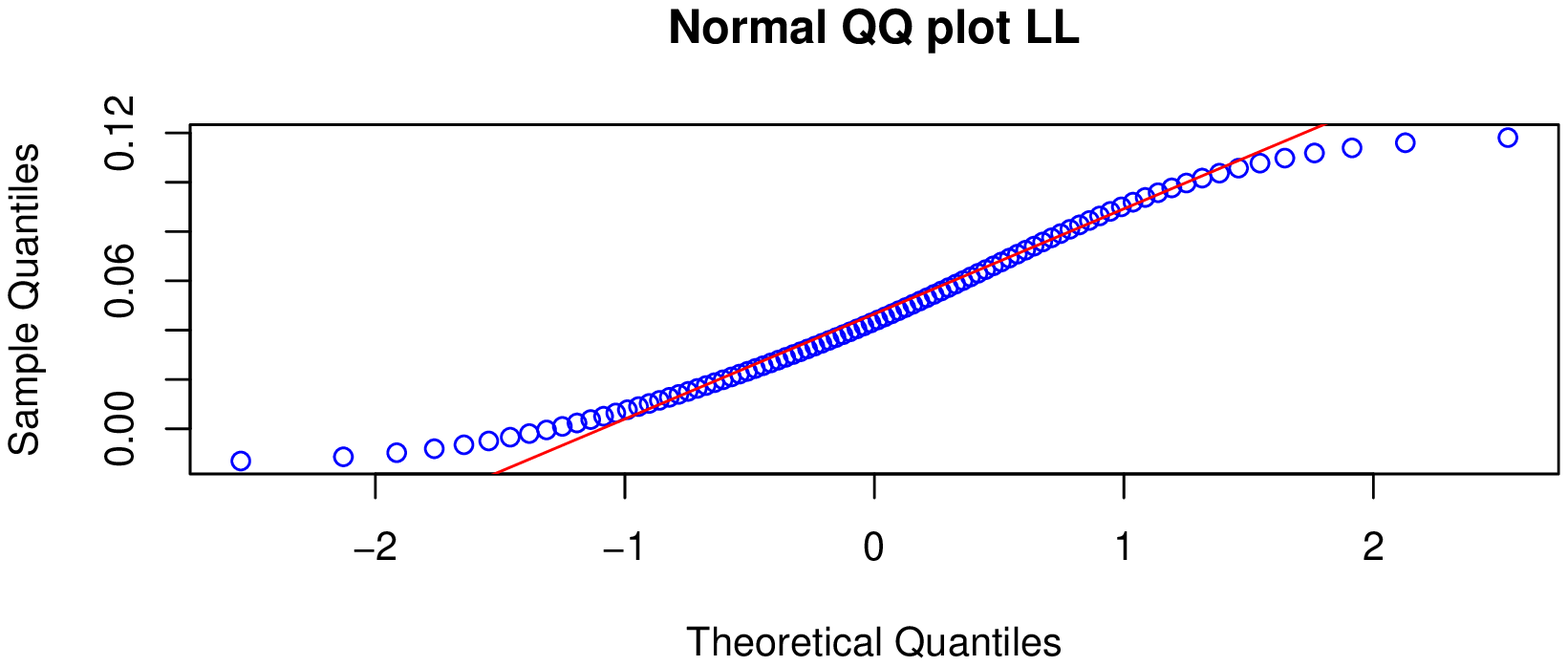}
\caption{ QQ plot of local linear estimator for $\mu(x) = - 10*x$
with $T = 10,~n = 1000,~\lambda = 2$ and jump size $Z_{n} \sim
\mathscr{N}(0,0.036^{2})$ }
\end{figure}

Next, we will assess the global performance between local linear
estimator and NW estimator via the Root of Mean Square Errors (RMSE)
\begin{equation}
\label{mc5} RMSE =
\sqrt{\frac{1}{m}\sum_{k=0}^{m}\left\{\hat{\mu}(x_{k}) -
\mu(x_{k})\right\}^{2}}, \end{equation} where $\hat{\mu}(x)$ is the
estimator of $\mu(x)$ and $\{x_{k}\}_{1}^{m}$ are chosen uniformly
to cover the range of sample path of $X_{t}.$ Tables 2, 3 and 4
report the results on RMSE-LL and RMSE-NW for the drift function
$\mu(x)$ with different types of time spans, sampling numbers, jump
intensities and jump sizes over 100 replicates.

We can notice that the local linear estimator performs a little
better (approximately reduced by half) than the NW estimator in
terms of the RMSE under different types of time spans, sampling
numbers, jump intensities and jump sizes. This fact confirms that
local linear estimator possesses the property of bias correction.
From Table 2, we can get the other two findings. Firstly, for the
same time interval $T$, as the sample sizes $n$ tends larger, the
performances of these estimators improved due to more information
used for estimators. Secondly, for the same sample sizes $n$, as the
time interval $T$ expands larger, the performances of these
estimators get worse due to the fact that more jumps happened in
larger time interval $T.$ From Table 3 and 4, for the same jump size
or the same jump arrival intensity, as the sample sizes $n$ tends
larger, the performances of the estimators for $\mu(x)$ are also
improved due to the fact that more jump information for estimation
procedure is collected as as $\Delta_{n} \rightarrow 0.$ However,
for the same sample sizes $n$, as the amplitude or frequency of jump
becomes larger, the RMSE of the estimators gradually becoming
larger.

\begin{table}[!htb]
  \centering
  \caption{Simulation results on RMSE-NW and RMSE-LL
for three lengths of time interval (T) and three sample sizes for
$\mu(x) = -10 x$ with jump size $Z_{n} \sim
\mathscr{N}(0,0.036^{2})$ and jump intensity $\lambda = 2$ over 100
replicates.}
    \begin{tabular}{ccccc}
    \hline
    Time Span & Estimators & $n=500$ & $n=1000$ & $n=2500$ \\
    \hline
    {$T = 10$} & RMSE-NW & 0.2566  & 0.1403  & 0.1002  \\
          & RMSE-LL & 0.1532  & 0.0813  & 0.0560  \\
    {$T = 20$} & RMSE-NW & 0.3734  & 0.2246  & 0.1334  \\
          & RMSE-LL & 0.1033  & 0.0970  & 0.0947  \\
    {$T = 50$} & RMSE-NW & 0.9885  & 0.3907  & 0.1522  \\
          & RMSE-LL & 0.2114  & 0.1492  & 0.0869  \\
    \hline
    \end{tabular}%
  \label{tab:addlabel}%
\end{table}%

\begin{table*}[!htb]
\centering \caption{Simulation results on RMSE-NW and RMSE-LL for
three types of jump intensity $\lambda$ and three sample sizes for
$\mu(x) = -10 x$ with $T = 10$ and jump size $Z_{n} \sim
\mathscr{N}(0, 0.036^{2})$ over 100 replicates.}
\begin{tabular}{ccccc}
\hline Jump Intensity & Estimators & $n=500$ & $n=1000$ & $n=2500$\\
\hline
$\lambda = 1$ & RMSE-NW & 0.2297  & 0.1117  & 0.0629  \\
          & RMSE-LL & 0.0896  & 0.0892  & 0.0480  \\
$\lambda = 2$ & RMSE-NW & 0.2566  & 0.1403  & 0.1002  \\
          & RMSE-LL & 0.1532  & 0.0813  & 0.0560  \\
$\lambda = 5$ & RMSE-NW & 0.3174  & 0.1763  & 0.1305  \\
          & RMSE-LL & 0.1603  & 0.0950  & 0.0736  \\
\hline
\end{tabular}
\end{table*}

\begin{table*}[!htb]
\centering \caption{Simulation results on RMSE-NW and RMSE-LL for
three types of jump size $Z_{n}$ and three sample sizes for $\mu(x)
= -10 x$ with $T = 10$ and jump intensity $\lambda = 2$ over 100
replicates.}
\begin{tabular}{ccccc}
\hline Jump Size & Estimators & $n=500$ & $n=1000$ & $n=2500$\\
\hline $Z_{n} \sim
\mathscr{N}(0, 0.036^{2})$ & RMSE-NW & 0.2566  & 0.1403  & 0.1002  \\
          & RMSE-LL & 0.1532  & 0.0813  & 0.0560  \\
$Z_{n} \sim
\mathscr{N}(0, 1)$ & RMSE-NW & 0.4532  & 0.2031  & 0.1690  \\
          & RMSE-LL & 0.2961  & 0.0931  & 0.0797  \\
$Z_{n} \sim
Cauchy(0, 1)$ & RMSE-NW & 0.7983  & 0.2156  & 0.1830  \\
          & RMSE-LL & 0.4083  & 0.1242  & 0.1538  \\
\hline
\end{tabular}
\end{table*}

{\textbf{Example 2 (Infinite Activity case).}} The infinite activity
jump $J_{t}$ in second-order diffusion model (\ref{mc1}) is a
Variance Gamma (VG) jump process, that is $J_{t} = c G_{t} + \eta
W_{G_{t}}$ with $c = -0.2,~\eta = 0.2$, where $G_{t}$ is an
independent Gamma process subject to Gamma(t/b,b) with $b =
var(G_{1}) = 0.23$ as that in Madan \cite{mad}. As is known, VG
process is an infinite activity jump process with finite variation.

For the infinite activity case, from Figures 4, 5, 6 and Tables 5, 6
we can observe the similar finding as the finite activity case,
which confirms the smaller bias and the normality in Theorems
\ref{mr}. This also shows that the methodology proposed in this
paper is robust to the presence of infinite activity jumps.
\begin{figure}[!htb]
\centering
\includegraphics[width=1 \textwidth]{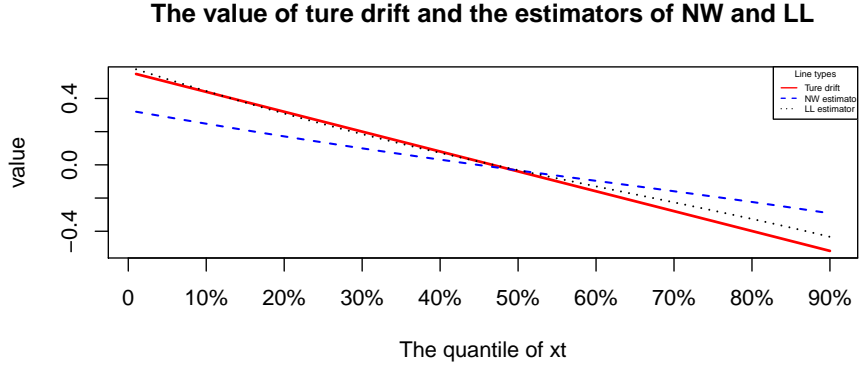}
\caption{ Various Nonparametric Estimators for $\mu(x) = - 10*x$
with $T = 10,~n = 1000$ and Variance Gamma jump process $J_{t}$ }
\end{figure}

\begin{table}[htbp]
  \centering
  \caption{The Biases of Nadaraya-Watson and Local Linear Estimators for $\mu(x) = -
  10*x$ at various quantile points of sample $X_t$ with $T = 10,~n = 1000$ and Variance Gamma jump process $J_{t}$}
    \begin{tabular}{cc|c|c|c|c|c|c|c|c}
    \toprule
    \multirow{2}[4]{*}{\textbf{Bias}} & \multicolumn{9}{c}{\textbf{Various Quantile Points of Sample Xt}} \\
\cmidrule{2-10}          & \multicolumn{1}{c}{\textbf{10\%}} & \multicolumn{1}{c}{\textbf{20\%}} & \multicolumn{1}{c}{\textbf{30\%}} & \multicolumn{1}{c}{\textbf{40\%}} & \multicolumn{1}{c}{\textbf{50\%}} & \multicolumn{1}{c}{\textbf{60\%}} & \multicolumn{1}{c}{\textbf{70\%}} & \multicolumn{1}{c}{\textbf{80\%}} & \textbf{90\%} \\
    \midrule
    \textbf{NW} & \multicolumn{1}{c}{-0.1702 } & \multicolumn{1}{c}{-0.1088 } & \multicolumn{1}{c}{-0.0564 } & \multicolumn{1}{c}{-0.0178 } & \multicolumn{1}{c}{-0.0110 } & \multicolumn{1}{c}{0.0418 } & \multicolumn{1}{c}{0.0732 } & \multicolumn{1}{c}{0.1045 } & 0.1590  \\
    \textbf{LL} & \multicolumn{1}{c}{-0.0708 } & \multicolumn{1}{c}{-0.0415 } & \multicolumn{1}{c}{-0.0190 } & \multicolumn{1}{c}{-0.0048 } & \multicolumn{1}{c}{-0.0039 } & \multicolumn{1}{c}{0.0109 } & \multicolumn{1}{c}{0.0152 } & \multicolumn{1}{c}{0.0159 } & 0.0069  \\
    \bottomrule
    \end{tabular}%
  \label{tab:addlabel}%
\end{table}%

\begin{table}[!htb]
  \centering
  \caption{Simulation results on RMSE-NW and RMSE-LL
for three lengths of time interval (T) and three sample sizes for
$\mu(x) = -10 x$ with Variance Gamma jump process $J_{t}$ over 100
replicates.}
    \begin{tabular}{ccccc}
    \toprule
    Time Span & Estimators & $n=500$ & $n=1000$ & $n=2500$ \\
    \midrule
    \multirow{2}[1]{*}{$T=10$} & RMSE-NW & 0.2493  & 0.1140  & 0.0962  \\
          & RMSE-LL & 0.1879  & 0.0933  & 0.0494  \\
    \multirow{2}[0]{*}{$T=20$} & RMSE-NW & 0.4365  & 0.2481  & 0.2239  \\
          & RMSE-LL & 0.2892  & 0.1564  & 0.1561  \\
    \multirow{2}[1]{*}{$T=50$} & RMSE-NW & 1.1106  & 0.6076  & 0.2544  \\
          & RMSE-LL & 0.3077  & 0.2624  & 0.1837  \\
    \bottomrule
    \end{tabular}%
  \label{tab:addlabel}%
\end{table}%

\begin{figure}[!htb]
\centering
\includegraphics[width=1 \textwidth]{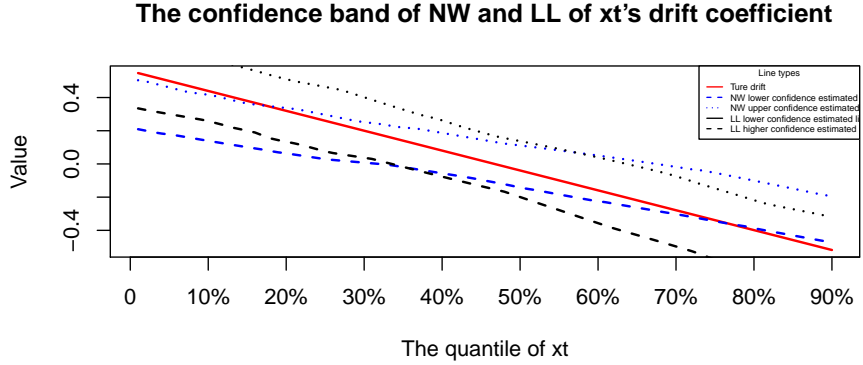}
\caption{ 95\% Monte Carlo Confidence Interval for $\mu(x) = - 10*x$
with $T = 10,~n = 1000$ and Variance Gamma jump process $J_{t}$ }
\end{figure}

\begin{figure}[!htb]
\centering
\includegraphics[width=1 \textwidth]{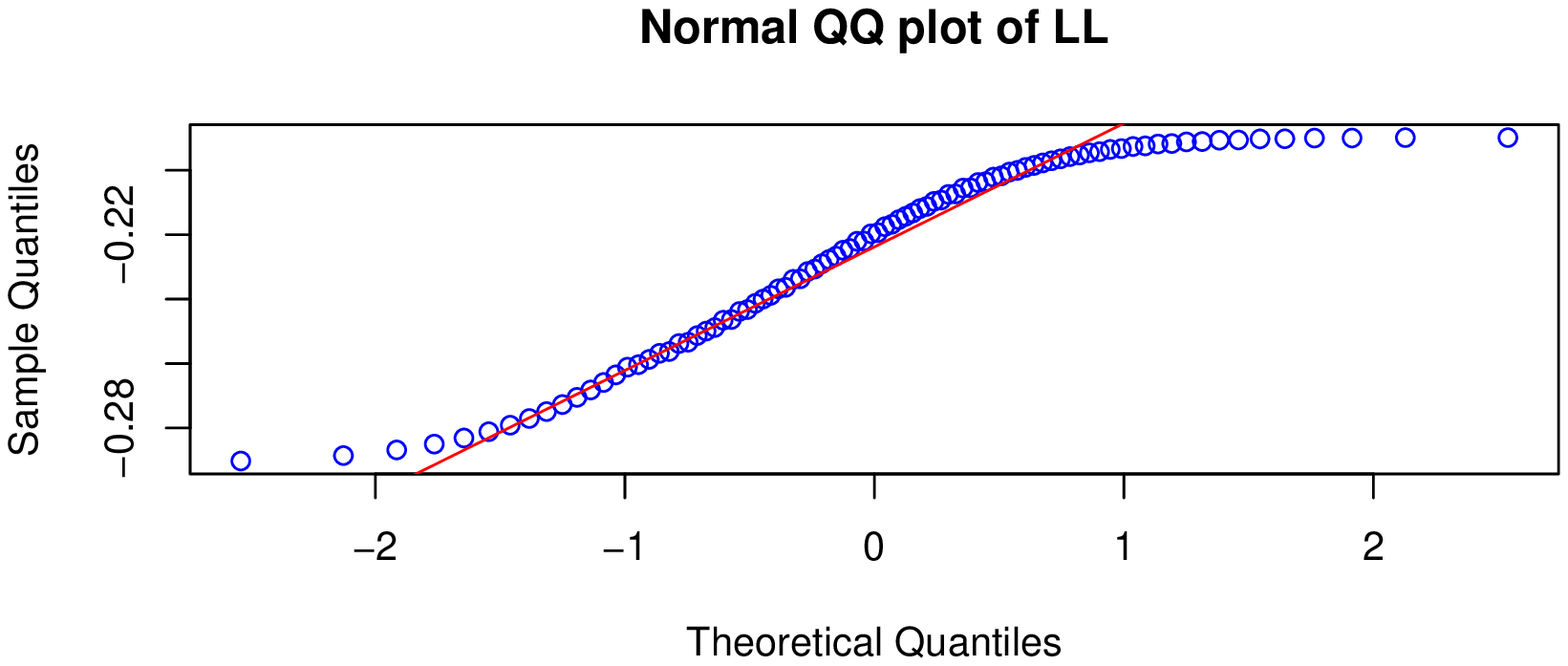}
\caption{ QQ plot of local linear estimator for $\mu(x) = - 10*x$
with $T = 10,~n = 1000$ and Variance Gamma jump process $J_{t}$ }
\end{figure}

\section{ \bf{Empirical Analysis }}
In this section, we apply the integrated diffusion with jump to
model the return of stock index from Shenzhen Stock Exchange in
China under five-minute high frequency data, that is $\Delta_{n} =
\frac{1}{48}$ (t = 1 meaning one day), and then apply the local
linear estimators to estimate the unknown coefficients in model
(\ref{om}).

We assume that
\begin{equation}
\label{mc6} \left\{
\begin{array}{ll}
d \log Y_{t} = X_{t}dt,\\
dX_{t} = \mu(X_{t-})dt + \sigma(X_{t-})dW_{t}+
\int_{\mathscr{E}}c(X_{t-} , z) r(\omega, dt, dz),
\end{array}
\right.\end{equation} where $\log Y_{t}$ is the log integrated
process for stock index or commodity price and $X_{t}$ is the latent
process for the log-returns. According to (\ref{yi}), we can get the
proxy of the latent process, that is the return of $\log Y_{t}$,
\begin{equation}
\label{mc7} \widetilde{X}_{i\Delta_n}=\frac{\log Y_{i\Delta_n} -
\log Y_{(i-1)\Delta_n}}{\Delta_n}.
\end{equation}

The plots of Shenzhen Composite Index and its proxy (\ref{mc7})
under five-minute frequency data from Jan 5, 2015 to Dec 31, 2015
are depicted in Figures 7 and 8. Figure 8 indicates the existence of
jumps.
\begin{figure}[!htb]
\centering
\includegraphics[width=1 \textwidth]{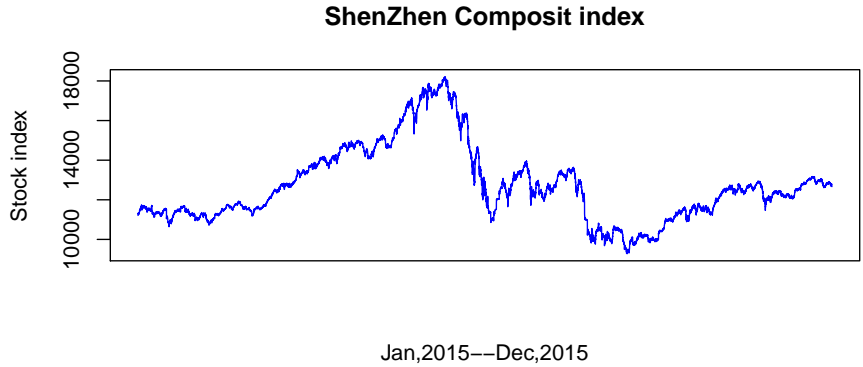}
\caption{ Time Series of Shenzhen Composite Index (2015) }
\end{figure}

\begin{figure}[!htb]
\centering
\includegraphics[width=1 \textwidth]{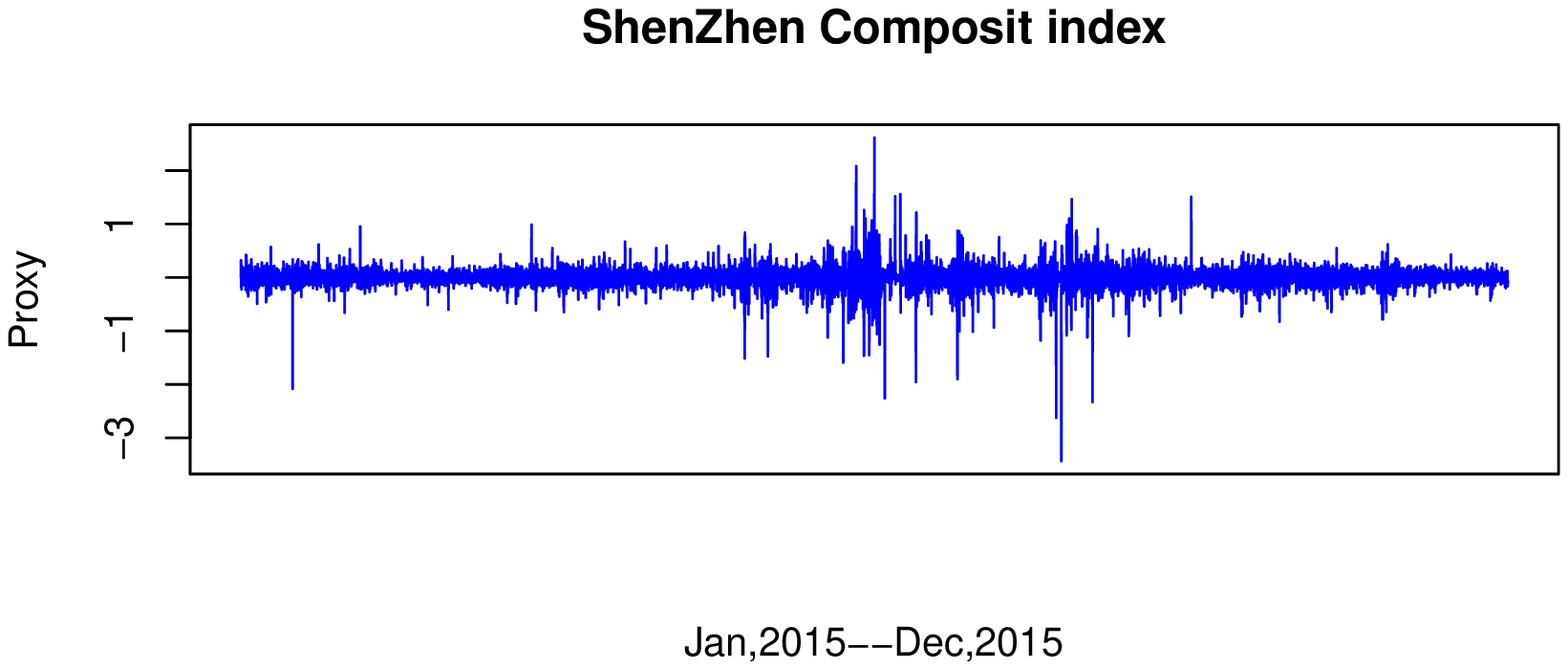}
\caption{ Proxy of Shenzhen Composite Index (2015) }
\end{figure}

Trough the Augmented Dickey-Fuller test statistic in Table 7, we can
easily observe that the null hypothesis of non-stationarity is
accepted for the logarithm of stock index $\log Y_{t}$ at the 5\%
significance level, but is rejected for the difference sequence
proxy $X_{t}$, which confirms the stationary Assumption \ref{a2} for
$X_{t}.$ Furthermore, based on the statistic proposed in
A\"it-Sahalia and Jacod \cite{aj}, the degree of activity jumps is
0.3487, which indicates the existence of infinite activity jumps for
$X_{t}$ and confirms the validity of model (\ref{mc6}) with infinite
activity jumps for Shenzhen Composite Index.

\begin{table*}[!htb]
\centering \caption{ Augmented Dickey-Fuller stationarity test for
model}
\begin{tabular}{crrrrrc}
\hline Data & \multicolumn{1}{c}{$~$} & \multicolumn{1}{c}{TestStat}
&
\multicolumn{1}{c}{CriticalValue} & \multicolumn{1}{c}{PValue}\\
\hline
  & $\log Y_{t}$ & -1.5298        & -2.8610        & 0.5046   \\
  & $X_{t}$      & -113.7941      & -2.8610        & $<0.01$    \\
\hline
\end{tabular}
\end{table*}

Here we use Gaussian kernels and the empirical bandwidth $h =
1.06*\hat{S}*T^{-1/5}$ for the estimation procedure. The local
linear estimators under (\ref{mc7}) for the unknown coefficients in
model (\ref{mc6}) are displayed in Figures 9 and 10. It is observed
that the linear shape with negative coefficient for drift estimator
in Figure 9 which reveals the economic phenomenon of mean reversion.
It is also shown the quadratic form with positive coefficient for
volatility estimator with a minimum at 0.07 in Figure 10, which
coincides with the economic phenomenon of volatility smile. The
shapes of estimated curves for these unknown coefficients coincide
with those in Nicolau \cite{n2}.

However there is an interesting phenomenon discovered in the
volatility line of return series: this line is asymmetric not like
the symmetric one in Nicolau \cite{n2} and the volatility has
different rates to positive return and negative return. The finding
coincides with the conclusion in most empirical analysis that
asymmetric features are depicted for volatility with respect to
positive perturbations (good news) and negative perturbations (bad
news). This is statistically due to the higher agglomeration effect
of jumps for Shenzhen stock index in 2015. As Johannes \cite{joh}
pointed out, jumps account for more than half of the conditional
variance. Furthermore, as Chen and Sun \cite{cs} concluded, the jump
behavior has a significant and asymmetric feedback effect in the
expected volatility, moreover, jumps exacerbate the degree of
asymmetric features of the volatility for stock index. Using the
statistics proposed in A\"it-Sahalia and Jacod \cite{aj}, we can
observe that at the points of the stock index return which are more
than 0.07, jumps happened at approximately 75.50\% points in them
and at the points of the return which are less than 0.07, jumps
happened at about 58.90\% points among them. In conclusion,
different frequencies of jumps at various points of stock index
return lead to the asymmetric volatility line.

What may explain this phenomenon economically is the special
situation of stock market in 2015. As we know, Shenzhen index in
China performed quite well during first half year of 2015, so more
shares changed hands arises for high yield, which leads to higher
slope for volatility at larger positive value of log-return
increments. However, Shenzhen index performed very unsatisfactorily
in the second half year. Facing the big transformation, investors
are still immersed in the past prosperity and believe that stocks
market crisis is temporary and the market can rebound after hitting
rock bottom. Hence less shares changed hands arises, which leads to
milder slope for volatility at lower value of log-return increments.
From Figure 11, we can observe that more and larger volume happened
at positive and large return rate. As a result, investors have
different responses to negative and positive returns, which leads to
the asymmetric volatility line.

\begin{figure}[!htb]
\centering
\includegraphics[width=1 \textwidth]{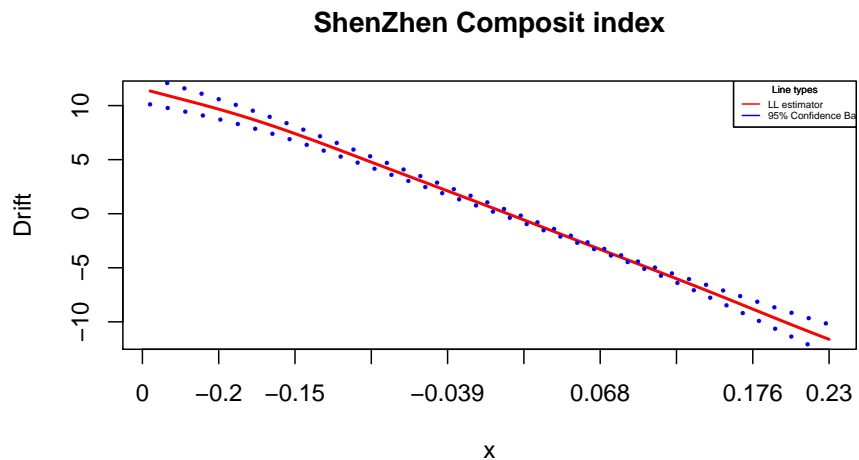}
\caption{ local linear estimator and its 95\% confidence bands of
the drift coefficient }
\end{figure}

\begin{figure}[!htb]
\centering
\includegraphics[width=1 \textwidth]{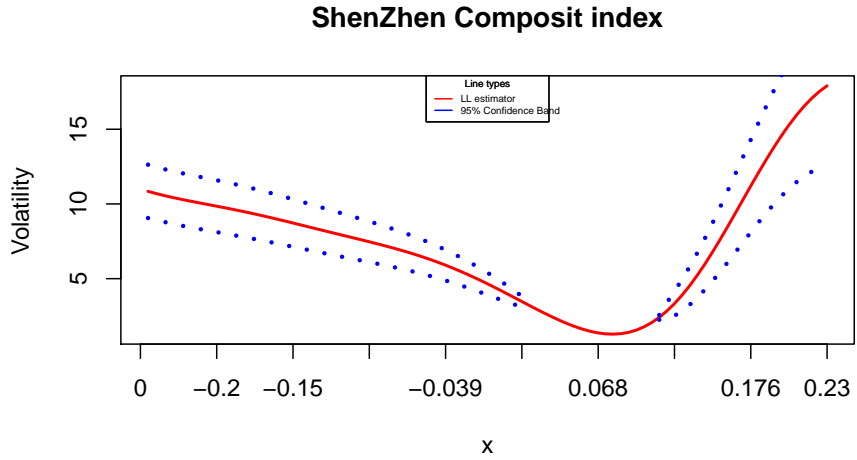}
\caption{ local linear estimator and its 95\% confidence bands of
the volatility coefficient }
\end{figure}

\begin{figure}[!htb]
\centering
\includegraphics[width=1 \textwidth]{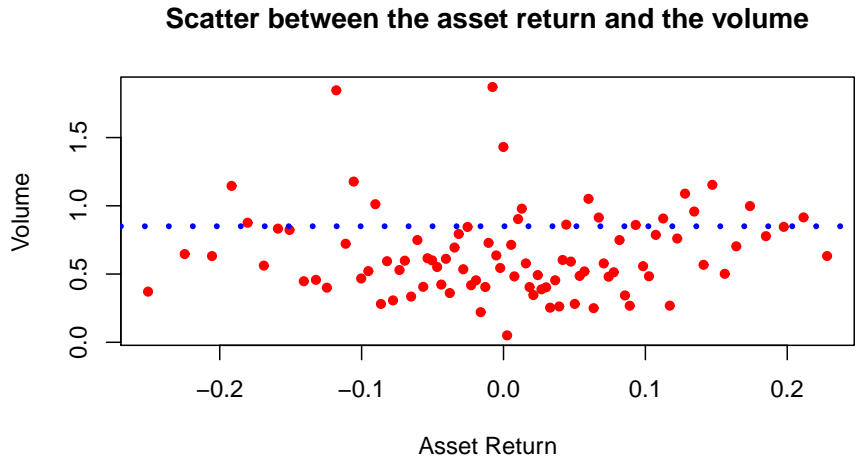}
\caption{ Scatter between the asset return rate and volume }
\end{figure}

\section{ Appendix }
In this section, we first present some technical lemmas and the
proofs for the main theorems.

\subsection{ Some Technical Lemmas with Proofs}

\begin{lemma}
\label{l1} (Shimizu and Yoshida \cite{syo}) Let $Z$ be a $d$-dimensional
solution-process to the stochastic differential equation $$Z_{t} =
Z_{0} + \int_{0}^{t}\mu(Z_{s-})ds + \int_{0}^{t}\sigma(Z_{s-})dW_{s}
+ \int_{0}^{t}\int_{\mathscr{E}}c(Z_{s-} , z)r(\omega, dt, dz),$$
where $Z_{0}$ is a random variable, $\mathscr{E} = \mathbb{R}^{d}
\setminus\{0\}$, $\mu(x), c(x , z)$ are $d$-dimensional vectors
defined on $\mathbb{R}^{d}, \mathbb{R}^{d}\times\mathscr{E}$
respectively, $\sigma(x)$ is a $d \times d$ diagnonal matrix defined
on $\mathbb{R}^{d}$, and $W_{t}$ is a $d$-dimensional vector of
independent Brownian motions.

Let $g$ be a $C^{2(l+1)}$-class function whose derivatives up to
2$(l+1)$th are of polynomial growth. Assume that the coefficient
$\mu(x), \sigma(x),$ and $c(x ,z)$ are $C^{2l}$-class function whose
derivatives with respective to $x$ up to 2$l$th are of polynomial
growth. Under Assumption 6, the following expansion holds
\begin{equation}
\label{3.1 } E[g(Z_{t})|\mathscr{F}_{s}] =
\sum_{j=0}^{l}L^{j}g(Z_{s})\frac{\Delta_{n}^{j}}{j!} + R,
\end{equation}for $t > s$ and $\Delta_{n} = t - s$, where  $R =
\int_{0}^{\Delta_{n}}\int_{0}^{u_{1}} \dots \int_{0}^{u_{l}}
E[L^{l+1}g(Z_{s+u_{l+1}})|\mathscr{F}_{s}]du_{1} \dots du_{l+1}$ is
a stochastic function of order $\Delta_{n}^{l+1}, Lg(x) =
\partial_{x}^{\ast}g(x)\mu(x) + \frac{1}{2}
tr[\partial_{x}^{2}g(x)\sigma(x)\sigma^{\ast}(x)] +
\int_{\mathscr{E}}\{g(x + c(x , z)) - g(x) -
\partial_{x}^{\ast}g(x)c(x , z)\}f(z)dz.$
\end{lemma}

\begin{remark}
\label{lr1} Consider a particularly important model:
$$\left\{
\begin{array}{ll}dY_{t} = X_{t-}dt,\\
dX_{t} = \mu(X_{t-})dt + \sigma(X_{t-})dW_{t} +
\int_{\mathscr{E}}c(X_{t-} , z)r(w , dt , dz).\end{array} \right.
$$
\noindent As $d$ = 2, we have
\begin{equation}
\label{do}
\begin{array}{ll}
Lg(x , y) = x(\partial g/\partial y) + \mu(x)(\partial g/\partial x)
+ \frac{1}{2}\sigma^{2}(x)(\partial^{2} g/\partial
x^{2})\\
~~~~~~~~~~~~~+ \int_{\mathscr{E}}\{g(x + c(x , z) , y) - g(x , y) -
\frac{\partial g}{\partial x} \cdot c(x , z) \}f(z)dz.
\end{array}
\end{equation}

Based on operator (\ref{do}), one can obtain the equations such as
(\ref{ice1}) and (\ref{ice2}), one can refer to Appendix A in Song,
Lin and Wang \cite{slw} for detailed calculations.
\end{remark}

\begin{lemma}
\label{l2} (Jacod \cite{ja12}) A sequence of $\mathbb{R}-$valued
variables $\{\zeta_{n, i}: i \geq 1\}$ defined on the filtered
probability space $(\Omega, \mathcal {F}, (\mathcal {F})_{t \geq 0},
\mathbb{P})$ is $\mathcal {F}_{i \Delta_{n}}-$measurable for all $n,
i.$ Assume there exists a continuous adapted $\mathbb{R}-$valued
process of finite variation $B_{t}$ and a continuous adapted and
increasing process $C_{t}$, for any $t
> 0,$ we have
\begin{equation}
\label{def4} \sup_{0 \leq s \leq t}\big|\sum_{i =
1}^{[s/\Delta_{n}]}\mathbb{E}\big[\zeta_{n, i} | \mathcal {F}_{(i -
1)\Delta_{n}}\big] - B_{s}\big| \stackrel{P} \longrightarrow 0,
\end{equation}
\begin{equation}
\label{def5} ~~~~~~~~~~~~~~~~~~~~~~~~~\sum_{i =
1}^{[t/\Delta_{n}]}\big(\mathbb{E}\big[\zeta_{n, i}^{2} | \mathcal
{F}_{(i - 1)\Delta_{n}}\big] - \mathbb{E}^{2}\big[\zeta_{n, i} |
\mathcal {F}_{(i - 1)\Delta_{n}}\big]\big) - C_{t} \stackrel{P}
\longrightarrow 0,
\end{equation}
\begin{equation}
\label{def6} \sum_{i = 1}^{[t/\Delta_{n}]}\mathbb{E}\big[\zeta_{n,
i}^{4} | \mathcal {F}_{(i - 1)\Delta_{n}}\big] \stackrel{P}
\longrightarrow 0.
\end{equation}
Then the processes $$\sum_{i = 1}^{[t/\Delta_{n}]}\zeta_{n, i}
\Rightarrow B_{t} + M_{t},$$ where $M_{t}$ is a continuous process
defined on the filtered probability space $\big({\Omega}, {P},
{\mathcal {F}}\big)$ and which, conditionally on the the
$\sigma-$filter $\mathcal {F}$, is a centered Gaussian
$\mathbb{R}-$valued process with $E \big[M_{t}^{2} | \mathcal
{F}\big] = C_{t}.$
\end{lemma}

\begin{remark}
Condition (\ref{def6}) is a conditional Lindeberg theorem or
Lyapounov's condition, whose aims is to ensure that the limiting
process is continuous. It is a particular case of Theorem
\uppercase\expandafter{\romannumeral9}.7.28 in Jacod and Shiryaev
\cite{js}.
\end{remark}

\begin{lemma}
\label{l3} (Song, Lin \cite{slw2}) Under Assumptions \ref{a1} -
\ref{a5}, we have
\begin{equation} \label{l3k}
\frac{1}{nh_{n}}\sum_{i=1}^{n}K\left(\frac{X_{i}-x}{h_{n}}\right)\left(\frac{X_{i}-x}{h_{n}}\right)^{k}~
\stackrel{\mathrm{a.s.}}{\longrightarrow}~
p(x)K_{1}^{k},\end{equation} where $K_{1}^{k} := \int u^{k}K(u) du.$
\end{lemma}
\begin{remark}
\label{lr3} We can obtain that
\begin{equation}
\frac{1}{nh_{n}^{3}}\sum_{j=1}^n K\left(\frac{X_{j-1} -
x}{h_n}\right)(X_{j-1} - x)^2
\stackrel{\mathrm{a.s.}}{\longrightarrow}~ p(x)K_{1}^{2},
\end{equation}
and
\begin{equation}
\frac{1}{nh_{n}^{2}}\sum_{j=1}^n K\left(\frac{X_{j-1} -
x}{h_n}\right)(X_{j-1} - x)
\stackrel{\mathrm{a.s.}}{\longrightarrow}~ p(x)K_{1}^{1}.
\end{equation}
\end{remark}

\begin{lemma}
\label{l4} Under Assumptions \ref{a1}- \ref{a5}, let
$$\mu_{n}^{\ast} (x)=\frac{\sum_{i=1}^n
w_{i-1}^{\ast}\left(\frac{X_i-X_{i-1}}{\Delta_{n}}\right)}{\sum_{i=1}^n
w_{i-1}^{\ast}}$$ and
$$ M_{n}^{\ast} (x)=\frac{\sum_{i=1}^n
w_{i-1}^{\ast}\frac{(X_i-X_{i-1})^2}{\Delta_{n}}}{\sum_{i=1}^n
w_{i-1}^{\ast}}.$$ where
$$w_i^ \ast =K\left(\frac{X_i - x}{h_n}\right)\left(\sum_{j=1}^n K\left(\frac{X_{j-1} - x}{h_n}\right)(X_{j-1} - x)^2-(X_i - x)
\sum_{j=1}^n K\left(\frac{X_{j-1} - x}{h_n}\right)(X_{j-1} -
x)\right)$$ then
$$\mu_{n}^{\ast}(x) \stackrel{p}{\rightarrow} \mu(x), ~~~ M_{n}^{\ast} (x) \stackrel{p}{\rightarrow} M(x).$$
Furthermore, if $h_{n} n \Delta_{n} \longrightarrow 0$ and $h_{n} =
O((n \Delta_{n})^{-1/5}),$ then
$$\sqrt{h_{n} n \Delta_n}\left(\mu_{n}^\ast (x) - \mu(x) - \frac{1}{2}h_{n}^{2}\mu^{''}(x)\frac{(K_{1}^{2})^{2} -
K_{1}^{3}K_{1}^{1}}{K_{1}^{2} - (K_{1}^{1})^{2}}\right)
\Longrightarrow N\left(0, V\frac{M(x)}{p(x)}\right),$$ and
$$\sqrt{h_{n} n \Delta_n}\left(M_{n}^{\ast} (x) - M(x) - \frac{1}{2}h_{n}^{2}M^{''}(x)\frac{(K_{1}^{2})^{2} -
K_{1}^{3}K_{1}^{1}}{K_{1}^{2} - (K_{1}^{1})^{2}}\right)
\Longrightarrow
N\left(0,V\frac{\int_{\mathscr{E}}c^4(x,y)\nu(dy)}{p(x)}\right),$$
where
$$V=\frac{(K_1^2)^2 K_2^0+ (K_1^1)^2 K_2^2 -2(K_1^1)(K_1^2)K_2^1}{[K_1^2-(K_1^1)^2]^2}.$$
\end{lemma}
\begin{remark}
\label{lr4} The main method to obtain the asymptotic properties for
estimator of model (\ref{om}) is to approximate the estimator for
(\ref{om}) by the similar estimator for univariate jump-diffusion
$X_{t}$ in probability. Based on the basic idea, we should first
know the asymptotic properties for the local linear estimator for
univariate jump-diffusion $X_{t}$ before the Theorems \ref{mr}
presented. Lemma \ref{l4} gives us the desired properties for local
linear estimator for univariate jump-diffusion $X_{t}.$
\end{remark}

\begin{proof}

As for the finite activity jumps case of $X_{t},$ Hanif \cite{hm1},
Chen and Zhang \cite{cz} considered the weak consistency and
asymptotic normality for the estimators $\mu_{n}^{\ast} (x)$ and
$M_{n}^{\ast} (x).$ Song, Lin \cite{slw2} discussed the weak
consistency and asymptotic normality of $\mu_{n}^{\ast} (x)$ for the
infinite activity jumps case of $X_{t}.$ They also derived the weak
consistency of $M_{n}^{\ast} (x)$ for the asymptotic variance of the
asymptotic distribution of $\mu_{n}^{\ast} (x).$ However, they did
not study the asymptotic distribution of $M_{n}^{\ast} (x)$ for the
infinite activity jumps case which we will give a detailed proof for
based on the Lemma \ref{l3}.

Using It\^{o} formula to the jump-diffusion setting shown in Protter
\cite{pr}, we can write
\begin{eqnarray*}
& ~ & (X_{(i+1)\Delta_{n}} - X_{i\Delta_{n}})^{2}\\
 & = & 2\int_{i\Delta_{n}}^{(i+1)\Delta_{n}}(X_{s-} -
 X_{i\Delta_{n}})\mu(X_{s-})ds + 2\int_{i\Delta_{n}}^{(i+1)\Delta_{n}}(X_{s-} -
 X_{i\Delta_{n}})\sigma(X_{s-})dW_{s}\\
 & ~ & + 2\int_{i\Delta_{n}}^{(i+1)\Delta_{n}}(X_{s-} -
 X_{i\Delta_{n}})\xi(X_{s-})dJ_{s} +
 \int_{i\Delta_{n}}^{(i+1)\Delta_{n}}(\sigma^{2}(X_{s-}) +
 \xi^{2}(X_{s-})Var(J(1)))ds\\
 & ~ & +
 \int_{i\Delta_{n}}^{(i+1)\Delta_{n}}\xi^{2}(X_{s-})\int_{\mathbb{R}}y^{2}\bar{\mu}(dy,
 ds).
\end{eqnarray*}
Hence we can derive that
\begin{eqnarray*}
& ~ & M_{n}^{\ast} (x) - M(x)\\ & = &
\frac{\frac{1}{n^{2}h_{n}^{4}}\sum_{i=1}^n
w_{i-1}^{\ast}\left((X_i-X_{i-1})^2 -
M(x) \cdot \Delta_{n}\right)}{\frac{\Delta_{n}}{n^{2}h_{n}^{4}}\sum_{i=1}^n w_{i-1}^{\ast}}\\
& = & \frac{1}{\frac{\Delta_{n}}{n^{2}h_{n}^{4}}\sum_{i=1}^n
w_{i-1}^{\ast}}
* \frac{1}{n^{2}h_{n}^{4}}\sum_{i=1}^n w_{i-1}^{\ast} \Bigg\{\left[\int_{(i-1)\Delta_{n}}^{i\Delta_{n}}(\sigma^{2}(X_{s-}) +
 \xi^{2}(X_{s-})Var(J(1)))ds - M(x) \cdot \Delta_{n}\right]\\
& ~ & + \left[2\int_{(i-1)\Delta_{n}}^{i\Delta_{n}}(X_{s-} -
 X_{i\Delta_{n}})\mu(X_{s-})ds\right] + \left[2\int_{(i-1)\Delta_{n}}^{i\Delta_{n}}(X_{s-} -
 X_{i\Delta_{n}})\sigma(X_{s-})dW_{s}\right]\\
 & ~ & + \left[2\int_{(i-1)\Delta_{n}}^{i\Delta_{n}}(X_{s-} -
 X_{i\Delta_{n}})\xi(X_{s-})dJ_{s}\right] + \left[\int_{(i-1)\Delta_{n}}^{i\Delta_{n}}\xi^{2}(X_{s-})\int_{\mathbb{R}}y^{2}\bar{\mu}(dy,
 ds)\right]\Bigg\}\\
& =: & b_{1_{n}} + A_{2_{n}} + A_{3_{n}} + A_{4_{n}} + A_{5_{n}}.
\end{eqnarray*}
We can express $b_{1_{n}},$ which dominates the bias for the
estimator of $M_{n}^{\ast} (x),$ as
\begin{eqnarray*}
b_{1_{n}} & = &
\frac{1}{\frac{\Delta_{n}}{n^{2}h_{n}^{4}}\sum_{i=1}^n
w_{i-1}^{\ast}} * \frac{1}{n^{2}h_{n}^{4}} \sum_{i=1}^n
w_{i-1}^{\ast}
\int_{(i-1)\Delta_{n}}^{i\Delta_{n}}(M(X_{s-}) - M(x))ds\\
& = & \frac{1}{\frac{\Delta_{n}}{n^{2}h_{n}^{4}}\sum_{i=1}^n
w_{i-1}^{\ast}} * \frac{1}{n^{2}h_{n}^{4}} \sum_{i=1}^n
w_{i-1}^{\ast} \int_{(i-1)\Delta_{n}}^{i\Delta_{n}}(M(X_{s-}) -
M(X_{(i-1)\Delta_{n}}))ds\\ & ~ & +
\frac{1}{\frac{1}{n^{2}h_{n}^{4}}\sum_{i=1}^n w_{i-1}^{\ast}} *
\frac{1}{n^{2}h_{n}^{4}} \sum_{i=1}^n
w_{i-1}^{\ast} (M(X_{(i-1)\Delta_{n}}) - M(x))\\
& = & \frac{1}{\frac{\Delta_{n}}{n^{2}h_{n}^{4}}\sum_{i=1}^n
w_{i-1}^{\ast}} * \frac{1}{n^{2}h_{n}^{4}} \sum_{i=1}^n
w_{i-1}^{\ast} \int_{(i-1)\Delta_{n}}^{i\Delta_{n}}(X_{s-} -
X_{(i-1)\Delta_{n}}) M^{'}(\tilde{X}_{s-}) ds\\
 & ~ & +
\frac{1}{\frac{1}{n^{2}h_{n}^{4}}\sum_{i=1}^n w_{i-1}^{\ast}} *
\frac{1}{n^{2}h_{n}^{4}} \sum_{i=1}^n w_{i-1}^{\ast}
[M^{'}(x)(X_{(i-1)\Delta_{n}} - x) +
\frac{1}{2}M^{''}(\hat{X}_{s-})(X_{(i-1)\Delta_{n}} - x)^{2}],
\end{eqnarray*}
where $\tilde{X}_{s-}$ lies between $X_{s-}$ and
$X_{(i-1)\Delta_{n}},$ $\hat{X}_{s-}$ lies between
$X_{(i-1)\Delta_{n}}$ and $x.$

Under a simple calculus, we have
\begin{equation}
\sum_{i=1}^n w_{i-1}^{\ast} (X_{(i-1)\Delta_{n}} - x) \equiv 0,
\end{equation}
so we can obtain that
\begin{eqnarray*}
b_{1_{n}} & = &
\frac{1}{\frac{\Delta_{n}}{n^{2}h_{n}^{4}}\sum_{i=1}^n
w_{i-1}^{\ast}} * \frac{1}{n^{2}h_{n}^{4}} \sum_{i=1}^n
w_{i-1}^{\ast} \int_{(i-1)\Delta_{n}}^{i\Delta_{n}}(X_{s-} -
X_{(i-1)\Delta_{n}}) M^{'}(\tilde{X}_{s-}) ds\\
 & ~ & +
\frac{1}{2} \frac{1}{\frac{1}{n^{2}h_{n}^{4}}\sum_{i=1}^n
w_{i-1}^{\ast}} * \frac{1}{n^{2}h_{n}^{4}} \sum_{i=1}^n
w_{i-1}^{\ast} M^{''}(\hat{X}_{s-})(X_{(i-1)\Delta_{n}} - x)^{2}\\
& =: & b_{11_{n}} + b_{12_{n}}.
\end{eqnarray*}
Now, define
$$\delta_{n, T}  = \max_{i \leq n} \sup_{i\Delta_{n} \leq s \leq (i+1)\Delta_{n}}\left|X_{s-} - X_{i\Delta_{n}}\right|.$$
Using the similar procedure as Bandi and Nguyen \cite{bn}, one can
have
\begin{equation}
\label{holder} \varlimsup_{n \rightarrow \infty} {\frac{\delta_{n,
T}}{\left(\Delta_{n} \log(1/\Delta_{n})\right)^{1/2}}} = C_{1}
~~~~~~ a.s.
\end{equation}
for some constant $C_{1},$ which implies that $\delta_{n, T} =
o_{a.s.}(1).$

It follows from the similar procedure of Lemma \ref{l3} as that in
Song, Lin \cite{slw2} for the numerator of $b_{12_{n}}$ with the
equation (\ref{holder})
\begin{equation}
\frac{1}{h_{n}^{2}} b^{Num}_{12_{n}}
\stackrel{\mathrm{a.s.}}{\longrightarrow} \frac{1}{2}
p^{2}(x)M^{''}(x)\left[(K_{1}^{2})^{2} - K_{1}^{3}K_{1}^{1}\right].
\end{equation}
With the asymptotic equations in Remark \ref{lr3}, we can easily get
for the denominator of $b_{12_{n}}$
\begin{equation}
b^{Den}_{12_{n}} = \frac{1}{n^{2}h_{n}^{4}}\sum_{i=1}^n
w_{i-1}^{\ast} \stackrel{\mathrm{a.s.}}{\longrightarrow}
p^{2}(x)\left[K_{1}^{2} - (K_{1}^{1})^{2}\right].
\end{equation}
Using the property that if $\xi_{n} \stackrel{P}{\rightarrow} \xi$
and $\eta_{n} \rightarrow c$ where $c$ is a nonzero constant, then
$\xi_{n}/\eta_{n} \stackrel{P}{\rightarrow} \xi/c,$ we have
\begin{equation}
\frac{1}{h_{n}^{2}} b_{12_{n}} \stackrel{P}{\longrightarrow}
\frac{1}{2}M^{''}(x)\frac{(K_{1}^{2})^{2} -
K_{1}^{3}K_{1}^{1}}{K_{1}^{2} - (K_{1}^{1})^{2}}.
\end{equation}
$b_{11_{n}}$ and $A_{2_{n}}$ can be dealt with the similar proof
procedure, here we only consider $A_{2_{n}}.$ Based on the equation
(\ref{holder}), the asymptotic equations in Remark \ref{lr3} and
Assumption \ref{a5}, we have that
\begin{equation}
A_{2_{n}} = \left(\Delta_{n}
\log(1/\Delta_{n})\right)^{1/2}O_{a.s.}(1) = o_{p}(b_{12_{n}}).
\end{equation}
Note that
\begin{equation}
A_{3_{n}} = \left(\Delta_{n}
\log(1/\Delta_{n})\right)^{1/2}O_{p}(A_{5_{n}}) = o_{p}(A_{5_{n}}),
\end{equation}
and
\begin{equation}
A_{4_{n}} = \left(\Delta_{n}
\log(1/\Delta_{n})\right)^{1/2}O_{p}(A_{5_{n}}) = o_{p}(A_{5_{n}}),
\end{equation}
which implies that $A_{5_{n}}$ is responsible for the asymptotic
distribution.

Now we only focus on the examination of the term $A_{5_{n}}.$ Denote
that
\begin{equation}
w_i^{\dag} = K\left(\frac{X_i -
x}{h_n}\right)\left(K_{1}^{2}-\left(\frac{X_i -
x}{h_{n}}\right)K_{1}^{1}\right),
\end{equation}
and
\begin{equation}
\frac{\sqrt{h_{n} n \Delta_n}}{n h_{n} \Delta_{n}} \sum_{i=1}^n
w_{i-1}^{\dag}
\int_{(i-1)\Delta_{n}}^{i\Delta_{n}}\xi^{2}(X_{s-})\int_{\mathbb{R}}y^{2}\bar{\mu}(dy,
 ds) := \sum_{i=1}^n q_{i}.
\end{equation}
We can obtain that $\sum_{i = 1}^{n}q_{i} \stackrel{D}
\longrightarrow N(0, [(K_{1}^{2})^{2}K_{2}^{0} -
2K_{1}^{1}K_{1}^{2}K_{2}^{1} +
(K_{1}^{1})^{2}K_{2}^{2}]\xi^{4}(x)\int_{\mathbb{R}}y^{4}\nu(dy)p(x))$
by Lemma \ref{l2}, if the following conditions hold

$|S_{1}| = \big| \sum_{i=1}^{n}E_{i-1}[q_{i}] \big|
\stackrel{P}{\rightarrow}0;$ \medskip

$S_{2} = \sum_{i=1}^{n}\big(E_{i-1}[{q_{i}}^{2}] - E_{i}^{2}[q_{1,
i+1}]\big) \stackrel{P}{\longrightarrow} [(K_{1}^{2})^{2}K_{2}^{0} -
2K_{1}^{1}K_{1}^{2}K_{2}^{1} +
(K_{1}^{1})^{2}K_{2}^{2}]\xi^{4}(x)\int_{\mathbb{R}}y^{4}\nu(dy)p(x);$\medskip

$S_{3} = \sum_{i=1}^{n}E_{i-1}[{q_{i}}^{4}]
\stackrel{P}{\rightarrow}0;$ \medskip

\noindent where $E_{i-1}[~\cdot~] = E[~\cdot~| X_{(i-1)\Delta_{n}}]$
and $q_{i}$ is $\mathcal {F}_{i \Delta_{n}}-$measurable, which is
the $\sigma-$algebra generated by $\{X_{j\Delta_{n}},~j \leq i\}.$

For $S_{2},$ the asymptotic variance can analogously be derived as
the drift case in Funke \cite{fk} by replacing $K\left(\frac{X_i -
x}{h_n}\right)$ with $w_i^{\dag},$ $\xi(X_{s-})$ with
$\xi^{2}(X_{s-})$ and $dL_{t} = \int_{\mathbb{R}}y\bar{\mu}(dy,
 ds)$ with $\int_{\mathbb{R}}y^{2}\bar{\mu}(dy, ds).$

For $|S_{1}|,$ since
\begin{eqnarray*}
E_{i-1}[q_{i}] & = & \frac{1}{\sqrt{h_{n} n \Delta_n}}
w_{i-1}^{\dag}
E_{i-1}\left[\int_{(i-1)\Delta_{n}}^{i\Delta_{n}}\xi^{2}(X_{s-})\int_{\mathbb{R}}y^{2}\bar{\mu}(dy,
 ds)\right]\\
& \equiv & 0,
\end{eqnarray*}
we can get $|S_{1}| \equiv 0.$

For $S_{3},$ with BDG inequality we obtain that
\begin{eqnarray*}
& ~ & \sum_{i=1}^{n} E_{i-1}[{q_{i}}^{4}]\\
& = & \frac{1}{(n h_{n} \Delta_{n})^{2}} \sum_{i=1}^{n}
w_{i-1}^{\dag 4}
E_{i-1}\left[\left(\int_{(i-1)\Delta_{n}}^{i\Delta_{n}}\xi^{2}(X_{s-})\int_{\mathbb{R}}y^{2}\bar{\mu}(dy,
 ds)\right)^{4}\right]\\
& \leq & C \frac{n}{(n h_{n} \Delta_{n})^{2}} (\Delta_{n}^{2} +
\Delta_{n})\\
& = & C \frac{1}{n h^{2}_{n}} + C \frac{1}{n \Delta_{n} h^{2}_{n}}
\rightarrow 0.
\end{eqnarray*}
Furthermore, we can obtain
\begin{eqnarray*}
& ~ & \frac{\sqrt{h_{n} n \Delta_n}}{n^{2}h_{n}^{4} \Delta_{n}}
\sum_{i=1}^n w_{i-1}^{\ast}
\int_{(i-1)\Delta_{n}}^{i\Delta_{n}}\xi^{2}(X_{s-})\int_{\mathbb{R}}y^{2}\bar{\mu}(dy,
 ds)\\
& - & \frac{\sqrt{h_{n} n \Delta_n}}{n h_{n} \Delta_{n}}
\sum_{i=1}^n w_{i-1}^{\dag}
\int_{(i-1)\Delta_{n}}^{i\Delta_{n}}\xi^{2}(X_{s-})\int_{\mathbb{R}}y^{2}\bar{\mu}(dy,
 ds)\\
& = & \left(\frac{1}{nh_{n}^{3}}\sum_{j=1}^n K\left(\frac{X_{j-1} -
x}{h_n}\right)(X_{j-1} - x)^2 - K_{1}^{2}\right) \times \\
& ~ & \times \frac{1}{\sqrt{h_{n} n \Delta_n}}\sum_{i=1}^n
K\left(\frac{X_i - x}{h_n}\right)
\int_{(i-1)\Delta_{n}}^{i\Delta_{n}}\xi^{2}(X_{s-})\int_{\mathbb{R}}y^{2}\bar{\mu}(dy,
 ds)\\
& ~ & - \left(\frac{1}{nh_{n}^{2}}K\left(\frac{X_{j-1} -
x}{h_n}\right)(X_{j-1} - x) - K_{1}^{1}\right) \times \\
& ~ & \times \frac{1}{\sqrt{h_{n} n \Delta_n}} \sum_{i=1}^n
K\left(\frac{X_i - x}{h_n}\right) \left(\frac{X_i - x}{h_n}\right)
\int_{(i-1)\Delta_{n}}^{i\Delta_{n}}\xi^{2}(X_{s-})\int_{\mathbb{R}}y^{2}\bar{\mu}(dy,
 ds)\\
& =: & \left(\frac{1}{nh_{n}^{3}}\sum_{j=1}^n K\left(\frac{X_{j-1} -
x}{h_n}\right)(X_{j-1} - x)^2 - K_{1}^{2}\right) \sum_{i=1}^n
q_{i}^{'}\\
& ~ & - \left(\frac{1}{nh_{n}^{2}}K\left(\frac{X_{j-1} -
x}{h_n}\right)(X_{j-1} - x) - K_{1}^{1}\right) \sum_{i=1}^n
q_{i}^{''}.
\end{eqnarray*}
Under the same proof procedure as $\sum_{i=1}^n q_{i},$ we can prove
that $\sum_{i=1}^n q_{i}^{'} \stackrel{D} \longrightarrow N(0,
K_{2}^{0}\xi^{4}(x)\int_{\mathbb{R}}y^{4}\nu(dy)p(x))$ and
$\sum_{i=1}^n q_{i}^{''} \stackrel{D} \longrightarrow N(0,
K_{2}^{2}\xi^{4}(x)\int_{\mathbb{R}}y^{4}\nu(dy)p(x)),$ which
implies that $\sum_{i=1}^n q_{i}^{'} = O_{p}(1)$ and $\sum_{i=1}^n
q_{i}^{''} = O_{p}(1).$ Hence with the Remark \ref{lr3} we can
obtain
\begin{eqnarray*}
& ~ & \frac{\sqrt{h_{n} n \Delta_n}}{n^{2}h_{n}^{4} \Delta_{n}}
\sum_{i=1}^n w_{i-1}^{\ast}
\int_{(i-1)\Delta_{n}}^{i\Delta_{n}}\xi^{2}(X_{s-})\int_{\mathbb{R}}y^{2}\bar{\mu}(dy,
 ds)\\
& - & \frac{\sqrt{h_{n} n \Delta_n}}{n h_{n} \Delta_{n}}
\sum_{i=1}^n w_{i-1}^{\dag}
\int_{(i-1)\Delta_{n}}^{i\Delta_{n}}\xi^{2}(X_{s-})\int_{\mathbb{R}}y^{2}\bar{\mu}(dy,
 ds) = o_{p}(1).
\end{eqnarray*}
By the Slutsky's theorem, we can get that
\begin{eqnarray*}
& ~ & \frac{\sqrt{h_{n} n \Delta_n}}{n^{2}h_{n}^{4} \Delta_{n}}
\sum_{i=1}^n w_{i-1}^{\ast}
\int_{(i-1)\Delta_{n}}^{i\Delta_{n}}\xi^{2}(X_{s-})\int_{\mathbb{R}}y^{2}\bar{\mu}(dy,
 ds) \\
 &\stackrel{D}
\longrightarrow& N(0, [(K_{1}^{2})^{2}K_{2}^{0} -
2K_{1}^{1}K_{1}^{2}K_{2}^{1} +
(K_{1}^{1})^{2}K_{2}^{2}]\xi^{4}(x)\int_{\mathbb{R}}y^{4}\nu(dy)p(x))
\end{eqnarray*}

Using the property that if $\xi_{n} \stackrel{D}{\rightarrow} \xi$
and $\eta_{n} \stackrel{P}{\rightarrow} c$ where $c$ is a nonzero
constant, then $\xi_{n}/\eta_{n} \stackrel{P}{\rightarrow} \xi/c,$
we have proved the asymptotic distribution of $M_{n}^{\ast} (x)$
\begin{equation}
\sqrt{h_{n} n \Delta_n} \left(M_{n}^{\ast} (x) - M(x) -
\frac{1}{2}h_{n}^{2}M^{''}(x)\frac{(K_{1}^{2})^{2} -
K_{1}^{3}K_{1}^{1}}{K_{1}^{2} - (K_{1}^{1})^{2}}\right)
\Longrightarrow
N\left(0,V\frac{\xi^{4}(x)\int_{\mathbb{R}}y^{4}\nu(dy)}{p(x)}\right),
\end{equation}
where $V=\frac{(K_1^2)^2 K_2^0+ (K_1^1)^2 K_2^2
-2(K_1^1)(K_1^2)K_2^1}{[K_1^2-(K_1^1)^2]^2}.$ \end{proof}

\subsection{ The proof of Theorem \ref{mr} }
\begin{proof}

As for the finite activity jumps case of $X_{t},$ Chen and Zhang
\cite{cz} considered the weak consistency and asymptotic normality
of the estimators $\hat{\mu}_{n}(x)$ and $\hat{M}_{n}(x)$ for the
second-order jump-diffusion. Here we will give a detailed proof for
the second-order diffusion with infinite activity jumps. For
simplicity, here we only prove the results for $\hat{\mu}_n(x)$ and
one can take the similar proof procedure for $\hat{M}_n(x).$

\noindent {\textbf{Weak Consistency:}} By Lemma \ref{l4}, it
suffices to show that :
$$ \hat{\mu}_n(x)-\mu^\ast _n(x) \stackrel{p}{\rightarrow} 0.$$
Firstly, we prove that
\begin{equation}
\label{6.10} \frac{1}{n^2 h_{n}^{4}}\sum_{i=1}^n
w_{i-1}-\frac{1}{n^2 h_{n}^{4}}\sum_{i=1}^n w_{i-1}^\ast
\stackrel{p}{\rightarrow} 0.
\end{equation}
To this end, we should prove that
\begin{equation}
\label{6.11} \frac{1}{n
h_{n}}\sum_{i=1}^{n}K\left(\frac{\tilde{X}_{i-1} - x}{h_{n}}\right)
- \frac{1}{n h_{n}}\sum_{i=1}^{n}K\left(\frac{X_{i-1} -
x}{h_{n}}\right) \stackrel{p}{\rightarrow} 0,
\end{equation}
\begin{equation}
\label{6.12}
\frac{1}{nh_{n}^{2}}\sum_{i=1}^{n}K\left(\frac{\tilde{X}_{i-1} -
x}{h_{n}}\right)(\tilde{X}_{i} - x) -
\frac{1}{nh_{n}^{2}}\sum_{i=1}^{n}K\left(\frac{X_{i-1} -
x}{h_{n}}\right)(X_{i-1} - x) \stackrel{p}{\rightarrow} 0,
\end{equation}
\begin{equation}
\label{6.13}
\frac{1}{nh_{n}^{3}}\sum_{i=1}^{n}K\left(\frac{\tilde{X}_{i-1} -
x}{h_{n}}\right)(\tilde{X}_{i} - x)^2-
\frac{1}{nh_{n}^{3}}\sum_{i=1}^{n}K\left(\frac{X_{i-1} -
x}{h_{n}}\right)(X_{i-1} - x)^2 \stackrel{p}{\rightarrow} 0.
\end{equation}
For (\ref{6.11}), let $\varepsilon_{1 , n} = \frac{1}{n
h_{n}}\sum_{i=1}^{n}K\left(\frac{\tilde{X}_{i-1} - x}{h_{n}}\right)
- \frac{1}{n h_{n}}\sum_{i=1}^{n}K\left(\frac{X_{i-1} -
x}{h_{n}}\right).$

\noindent By the mean-value theorem, stationarity, {Assumptions
\ref{a3}, \ref{a5} and (\ref{holder})}, we obtain
\begin{eqnarray*}
E[|\varepsilon_{1 , n}|] & \leq & E\left[\frac{1}{nh_{n}} \sum
_{i=1}^{n}|K ^{'}(\xi_{n,i}) \frac{\widetilde{X}_{(i-1)\Delta_{n}} -
{X}_{(i-1)\Delta_{n}}}{h_{n}}|\right]\\ & = &
E\left[\frac{1}{h_{n}}| K ^{'}(\xi_{n,2})
\frac{\widetilde{X}_{\Delta_{n}} -
{X}_{\Delta_{n}}}{h_{n}}|\right]\\& \leq
&{\frac{\sqrt{\Delta_{n}\log(1/\Delta_{n})}}{h_{n}}E\left[\frac{1}{h_{n}}|K
^{'}(\xi_{n,2})|\right] \rightarrow 0},
\end{eqnarray*}
where $\xi_{n,2}$ = $\theta(\frac{x - X_{\Delta_{n}}}{h}) + (1 -
\theta)(\frac{x - \widetilde{X}_{\Delta_{n}}}{h}) ~0 \leq \theta
\leq 1.$ Hence, (\ref{6.11}) follows from Chebyshev's inequality.

For (\ref{6.12}) we should prove that
\begin{equation}
\label{6.15}
\delta_{1,n}=\frac{1}{nh_{n}^{2}}\sum_{i=1}^{n}K\left(\frac{\tilde{X}_{i-1}
- x}{h_{n}}\right)(\tilde{X}_{i} - x) -
\frac{1}{nh_{n}^{2}}\sum_{i=1}^{n}K\left(\frac{{X}_{i-1} -
x}{h_{n}}\right)(\tilde{X}_{i} - x) \stackrel{p}{\rightarrow} 0,
\end{equation}
and
\begin{equation}
\label{6.16}
\delta_{2,n}=\frac{1}{nh_{n}^{2}}\sum_{i=1}^{n}K\left(\frac{{X}_{i-1}
- x}{h_{n}}\right)(\tilde{X}_{i} - x) -
\frac{1}{nh_{n}^{2}}\sum_{i=1}^{n}K\left(\frac{{X}_{i-1} -
x}{h_{n}}\right)({X}_{i - 1} - x) \stackrel{p}{\rightarrow} 0.
\end{equation}
Under Assumptions \ref{a3} and \ref{a5}, we have
\begin{eqnarray*}
\big|E[\delta_{1,n}]\big| & = & \left|E\left[ \frac{1}{nh_{n}^{2}}
\sum_{i=1}^{n} \left\{K\left(\frac{\tilde{X}_{i-1} -
x}{h_{n}}\right) - K\left(\frac{{X}_{i-1}
- x}{h_{n}}\right)\right\}(\tilde{X}_{i} - x)\right]\right|\\
& = &
\frac{1}{h_{n}^{2}}\left|E\left[\left\{K\left(\frac{\tilde{X}_{i-1}
- x}{h_{n}}\right)-K\left(\frac{{X}_{i-1} - x}{h_{n}}\right)\right\}
E[\tilde{X}_{i} - x\big|\mathscr{F}_{i-1}]\right]\right|\\
& = & \frac{1}{h_{n}^{2}} \left|E[ K'(\xi_{n,i})
\frac{(\tilde{X}_{i-1} - X_{i - 1})}{h_{n}} (X_{i-1} - x + O_P(\Delta_n))]\right|\\
& \leq & \frac{\sqrt{\Delta_{n} \log (1/\Delta_{n})}}{h_{n}^{2}} E[
\frac{1}{h_{n}}|K'(\xi_{n,i})
(X_{i-1} - x + O_P(\Delta_n))|]\\
& \rightarrow & 0,
\end{eqnarray*}
by stationarity, the mean-value theorem and Remark \ref{lr1}. So
$E[\delta_{1,n}]\rightarrow 0.$

If we can prove $Var[\delta_{1,n}] \rightarrow 0$, then (\ref{6.15})
holds. Now we calculate $Var[\delta_{1,n}]$
\begin{eqnarray*}
Var[\delta_{1,n}] & = &
\frac{1}{nh_{n}}Var[\frac{1}{\sqrt{n}}\sum_{i=1}^{n}
\frac{1}{h_{n}^{3/2}} K^{'}(\xi_{n,i})\frac{\widetilde{X}_{i-1}
- X_{i-1}}{h_{n}}(\widetilde{X}_{i} - x)]\\
& =: & \frac{1}{nh_{n}}Var[\frac{1}{\sqrt{n}}\sum_{i=1}^{n}f_{i}].
\end{eqnarray*}
where $f_{i} := \frac{1}{h_{n}^{3/2}}
K^{'}(\xi_{n,i})\frac{\widetilde{X}_{i-1} -
X_{i-1}}{h_{n}}(\widetilde{X}_{i} - x).$

\noindent By Remark \ref{lr1} and Assumptions \ref{a3} and \ref{a5},
we get
\begin{eqnarray*}
E[f_{i}^{2}] & = & E\left[\frac{1}{h_{n}^{3}}
K^{'2}(\xi_{n,i})\frac{(\widetilde{X}_{i-1} -
X_{i-1})^{2}}{h^{2}_{n}}E[(\widetilde{X}_{i} - x)^{2} | \mathscr{F}_{i - 1}]\right]\\
& \leq & \frac{\Delta_{n} \log (1/\Delta_{n})}{h^{2}_{n}}
E\left[\frac{1}{h_{n}^{3}} K^{'2}(\xi_{n,i}) ((X_{i} - x)^{2} +
O_{p}(\Delta_{n}))\right]\\
& < & \infty
\end{eqnarray*}
We notice that $f_i$ is stationary under Assumption \ref{a2} and
$\rho$-mixing with the same size as process
$\{\tilde{X}_{i\Delta_n}; i=1,2,...\}$ and $\{X_{i\Delta_n};
i=1,2,...\}$. So from Lemma 10.1.c with p=q=2 in Lin and
Bai(\cite{lb}, p. 132), we have
\begin{eqnarray*}
\big|Var[\frac{1}{\sqrt{n}}\sum_{i=1}^n f_i]\big|& = &
\big|\frac{1}{n}[\sum_{i=1}^n Var(f_i)+2
\sum_{j=1}^{n-1}\sum_{i=j+1}^n (Ef_if_j-Ef_iEf_j)]\big|\\
& = & Var(f_i)+\frac{2}{n}\big| \sum_{j=1}^{n-1}\sum_{i=j+1}^n
(Ef_if_j-Ef_iEf_j)]\big|\\
& \leq & Var(f_i)+ \frac{2}{n}\sum_{j=1}^{n-1}\sum_{i=j+1}^n
\big|Ef_if_j-Ef_iEf_j\big]\\
& \leq &
Var(f_i)+\frac{8}{n}\sum_{j=1}^{n-1}\sum_{i=j+1}^{n}\rho((i-j)\Delta_{n})(Ef^{2}_{i})^{\frac{1}{2}}(Ef^{2}_{j})^{\frac{1}{2}}\\
& = &
Var(f_i)+\frac{8}{n}\sum_{j=1}^{n-1}\sum_{i=j+1}^{n}\rho((i-j)\Delta_{n})Ef^{2}_{i}
\end{eqnarray*}
We have proved $Ef_i^2 < \infty$ above, so the first part in the
last equality $Var(f_i) < \infty $. Moreover, under Assumption
\ref{a2}, we have $\sum_{j=i+1}^n \rho((j-i)\Delta_n)
=O(\frac{1}{\Delta^{\alpha}_n})$. So $Var(\delta_{1,n}) =
\frac{1}{nh_n\Delta^{\alpha}_n}\rightarrow 0$ as $nh_n\Delta^{1 +
\alpha}_n\rightarrow \infty.$

Similar to the proof of (\ref{6.15}), we prove (\ref{6.16}) by
verifying $E[\delta_{2,n}]\rightarrow 0$ and
$Var[\delta_{2,n}]\rightarrow 0$. From the stationarity, Remark
\ref{lr1} and Assumptions \ref{a3} and \ref{a5}, we have
\begin{eqnarray*}
E[\delta_{2,n}] & = & \frac{1}{h_{n}^{2}}E\left[K\left(\frac{X_{i-1} - x}{h_{n}}\right)E\big[\widetilde{X}_{i} - X_{i-1}|\mathscr{F}_{i-1}\big]\right]\\
& = & \frac{\Delta_{n}}{h_{n}} E\left[\frac{1}{h_{n}}
K\left(\frac{X_{i-1} - x}{h_{n}}\right) \mu({X}_{i-1})\right] +
O\left(\frac{\Delta_{n}^{2}}{h^{2}_{n}}\right) \rightarrow 0.
\end{eqnarray*}
and
\begin{eqnarray*}
Var[\delta_{2,n}] & = & \frac{\Delta_{n}}{n h^{3}_{n}}
Var\left[\frac{1}{\sqrt{nh_{n}}}\sum_{i=1}^{n}
K\left(\frac{X_{i-1} - x}{h_{n}}\right) \frac{1}{\sqrt{\Delta_{n}}}(\widetilde{X}_{i} - X_{i-1})\right]\\
& =: & \frac{\Delta_{n}}{n h^{3}_{n}}
Var[\frac{1}{\sqrt{n}}\sum_{i=1}^{n}g_{i}],
\end{eqnarray*}
where
\begin{eqnarray*}
E\big[g_{i}^{2}\big] & = &
E\left[\frac{1}{h_{n}} K^{2}\left(\frac{X_{i-1} - x}{h_{n}}\right) E\left[\frac{(\widetilde{X}_{i} - X_{i-1})^{2}}{\Delta_{n}} | \mathscr{F}_{i - 1}\right]\right]\\
& = & \frac{1}{3}E\left[\frac{1}{h_{n}} K^{2}\left(\frac{X_{i-1} -
x}{h_{n}}\right) \big(\sigma^{2}({X}_{i-1}) +
\int_{\mathscr{E}}c^{2}({X}_{i-1}, z)f(z)dz\big)\right]\\
& \approx & O(1) < \infty.
\end{eqnarray*}
by the Remark \ref{lr1} and Assumptions \ref{a1}, \ref{a3}. Hence,
$Var[\delta_{2,n}]=O(\frac{\Delta^{1-\alpha}_n}{n h_n})\rightarrow
0$ under Assumption \ref{a5}. The proof of (\ref{6.13}) is similar
to that of (\ref{6.12}), so we omit it.

\noindent Secondly, we prove
\begin{equation}
\label{6.17} \delta_n:=\frac{1}{n^2 h_{n}^{4}}
\sum_{i=1}^nw_{i-1}\frac{\widetilde{X}_{i+1}-\widetilde{X}_{i}}{\Delta_n}
-\frac{1}{n^2
h_{n}^{4}}\sum_{i=1}^nw^\ast_{i-1}\frac{X_{i}-X_{i-1}}{\Delta_n}
\stackrel{p}{\rightarrow} 0.
\end{equation} which suffice to prove
\begin{equation}
\label{6.18} \frac{1}{n^2
h_{n}^{4}}\sum_{i=1}^nw^\ast_{i-1}\frac{\widetilde{X}_{i+1}-\widetilde{X}_{i}}{\Delta_n}
-\frac{1}{n^2
h_{n}^{4}}\sum_{i=1}^nw^\ast_{i-1}\frac{X_{i}-X_{i-1}}{\Delta_n}
\stackrel{p}{\rightarrow} 0,
\end{equation} and
\begin{equation}
\label{6.19} \frac{1}{n^2
h_{n}^{4}}\sum_{i=1}^nw_{i-1}\frac{\widetilde{X}_{i+1}-\widetilde{X}_{i}}{\Delta_n}
-\frac{1}{n^2
h_{n}^{4}}\sum_{i=1}^nw^\ast_{i-1}\frac{\widetilde{X}_{i+1}-\widetilde{X}_{i}}{\Delta_n}
\stackrel{p}{\rightarrow} 0.
\end{equation}

\noindent For (\ref{6.18}), we need only prove
\begin{equation}
\label{6.20} \delta_{3,n}=\frac{1}{n h_{n}}\sum_{i=1}^n
K\left(\frac{X_{i-1} - x}{h_{n}}\right)
\Big[\frac{\widetilde{X}_{i+1}-\widetilde{X}_{i}}{\Delta_n}-\frac{X_{i}-X_{i-1}}{\Delta_n}\Big]
\stackrel{p}{\rightarrow}0,
\end{equation} and
\begin{equation}
\label{6.21} \delta_{4,n}=\frac{1}{n h_{n}^{2}}\sum_{i=1}^n
K\left(\frac{X_{i-1} - x}{h_{n}}\right)(X_{i-1} - x)
\Big[\frac{\widetilde{X}_{i+1}-\widetilde{X}_{i}}{\Delta_n}-\frac{X_{i}-X_{i-1}}{\Delta_n}\Big]
\stackrel{p}{\rightarrow}0.
\end{equation}

\noindent the proof of (\ref{6.20}) and (\ref{6.21}) are similar, so
we just prove (\ref{6.21}).

By Lemma \ref{l1}, we can get
\begin{eqnarray*}
E[\varepsilon_{1,n}] & := & E\Big[\Big(\frac{(\widetilde{X}_{i+1} -
\widetilde{X}_{i})}{\Delta_{n}} - \frac{({X}_{i} -
{X}_{i-1})}{\Delta_{n}}\Big)\Big|\mathscr{F}_{i-1}\Big]\\
& = & E\Big\{E\Big[\Big(\frac{(\widetilde{X}_{i+1} -
\widetilde{X}_{i})}{\Delta_{n}} - \frac{({X}_{i} -
{X}_{i-1})}{\Delta_{n}}\Big)\Big|\mathscr{F}_{i}\Big]\Big|\mathscr{F}_{i-1}\Big\}\\
& = & \frac{\Delta_n}{2}\Big\{\mu(X_{i-1})\mu^{'}(X_{i-1})+\frac{1}{2}\sigma^2(X_{i-1})\mu^{''}(X_{i-1})\\
&~& + \int_{\mathscr{E}}\big\{\mu(X_{i-1}+ c(X_{i-1} , z)) -
\mu(X_{i-1}) - \mu^{'}(X_{i-1})\cdot c(X_{i-1},z)\big\}f(z)dz
\Big\},
\end{eqnarray*}
so by stationarity and Assumption \ref{a3}, we have
\begin{eqnarray*}
E[\delta_{4,n}(x)] & = &
E\Big\{\frac{1}{h_{n}^{2}}E\Big[K\left(\frac{X_{i-1} -
x}{h_{n}}\right)(X_{i-1} - x)\Big(\frac{\widetilde{X}_{i+1} -
\widetilde{X}_{i}}{\Delta_{n}} - \frac{{X}_{i} -
{X}_{i-1}}{\Delta_{n}}\Big)|\mathscr{F}_{i-1}\Big]\Big\}\\
& = &
\frac{\Delta_{n}}{2}E\Big[\frac{1}{h_{n}^{2}}K\left(\frac{X_{i-1} -
x}{h_{n}}\right)(X_{i-1} - x)
\Big(\mu(X_{i-1})\mu^{'}(X_{i-1}) +\frac{1}{2}\sigma^{2}(X_{i-1})\mu^{''}(X_{i-1})\\
& ~ &+ \int_{\mathscr{E}}\big\{\mu(X_{i-1} + c(X_{i-1} , z)) -
\mu(X_{i-1}) - \mu^{'}(X_{i-1})\cdot
c(X_{i-1},z)\big\}f(z)dz\Big)\Big]\\
& = & O(\Delta_{n})
\end{eqnarray*}
and
\begin{eqnarray*}
&&Var[\delta_{4,n}(x)] \\
& = & \frac{1}{n h_{n}
\Delta_{n}}Var\Big[\frac{1}{\sqrt{n}h^{3/2}_{n}}\sum_{i=1}^{n}{K\left(\frac{X_{i-1}
- x}{h_{n}}\right)\sqrt{\Delta_{n}}(X_{i-1} -
x)\Big(\frac{(\widetilde{X}_{i+1} - \widetilde{X}_{i})}{\Delta_{n}}
- \frac{({X}_{i\Delta_{n}} -
{X}_{i-1})}{\Delta_{n}}\Big)}\Big]\\
& =: & \frac{1}{n h_{n}
\Delta_{n}}Var\Big[\frac{1}{\sqrt{n}}\sum_{i=1}^{n}g_{i}\Big].
\end{eqnarray*}
\noindent By the similar analysis as above, we can easily obtains
$Var[\delta_{4,n}(x)]\rightarrow 0$ under Assumption \ref{a2} if
$E[g_i^2]<\infty$. In fact by Assumptions \ref{a1}, \ref{a3} and
\ref{a4}, we have
\begin{eqnarray*}
E[g_{i}^{2}] & = & E\Big[\frac{1}{h_{n}^{3}}
K^{2}\left(\frac{X_{i-1} - x}{h_{n}}\right)\Delta_{n}(X_{i-1} -
x)^2\Big(\frac{(\widetilde{X}_{i+1} -
\widetilde{X}_{i})}{\Delta_{n}} - \frac{({X}_{i} -
{X}_{i-1})}{\Delta_{n}}\Big)^{2}\Big]\\
& = & E\Big\{\frac{1}{h_{n}^{3}} K^{2}\left(\frac{X_{i-1} -
x}{h_{n}}\right)(X_{i-1} - x)^2
E\Big[\Delta_n\Big(\frac{\tilde{X}_{i+1}-\tilde{X}_{i}}{\Delta_n}-\frac{X_i-X_{i-1}}{\Delta_n}\Big)^2|\mathscr{F}_{i-1}\Big]\Big\}\\
& = & E\Big\{\frac{1}{h_{n}^{3}} K^{2}\left(\frac{X_{i-1} -
x}{h_{n}}\right)(X_{i-1} - x)^2
\times \frac{2}{3}\Big[\sigma^2(X_{i-1})+\int_{\mathscr{E}}c^2(X_{i-1},z)f(z)dz+O_P(\Delta_n)\Big]\Big\}\\
& \approx & O(1) < \infty.
\end{eqnarray*}
The proof of (\ref{6.19}) is similar to that of (\ref{6.18}).
Combination (\ref{6.10}) and (\ref{6.17}), the relationship
$\mu^\ast _n(x)-\hat{\mu}_n(x) \stackrel{p}{\rightarrow} 0$ holds,
 so by Lemma \ref{l1} we have $\hat{\mu}_n(x) \stackrel{p}{\rightarrow} \mu(x).$

\noindent {\textbf{Asymptotic Normality:}} By Lemma \ref{l4}, we
have
$$ U_n^\ast(x) := \sqrt{n\Delta_{n} h}\big(\mu_{n}^{\ast}(x) -
\mu(x) - h^{2} B_{\mu(x)}\big) \stackrel{d}{\rightarrow} N\Big( 0 ,
V\frac{M(x)}{p(x)}\Big),$$ where $B_{\mu(x)} =
\frac{1}{2}\mu^{''}(x)\frac{K_{2}^{2} -
K_{1}^{3}K_{1}^{1}}{K_{1}^{2} - (K_{1}^{1})^{2}}.$

 So by the asymptotic
equivalence theorem, it suffices to prove that
$$\hat{U}_n(x)-U_n^\ast(x) = \sqrt{h_{n} n \Delta_n }(\hat{\mu}_{n} (x)-\mu_{n}^\ast (x))\stackrel{p}{\rightarrow} 0,$$
where $\hat{U}_n(x) := \sqrt{n\Delta_{n} h}\big(\hat{\mu}_{n}(x) -
\mu(x) - h^{2} B_{\mu(x)}\big).$

In fact, from the proof of weak convergence such as (\ref{6.10}) and
(\ref{6.17}), we know that
\begin{eqnarray*}
& ~ & \hat{U}_n(x)-U_n^\ast(x)\\
& = &  \sqrt{h_{n} n \Delta_n }(\hat{\mu}_{n} (x)-\mu_{n}^\ast (x))\\
& = &  \sqrt{h_{n} n \Delta_n }\left(\frac{\delta_n}{\frac{1}{n^2
h_{n}^{4}}\sum_{i=1}^{n}w^\ast_{i-1}} +\frac{\sum_{i=1}^n w_{i-1}
\left(\frac{\widetilde{X}_{i+1}-\widetilde{X}_i}{\Delta_n}\right)}{\sum_{i=1}^{n}w_{i-1}}
-\frac{\sum_{i=1}^n w_{i-1} \left(\frac{\widetilde{X}_{i+1}-\widetilde{X}_i}{\Delta_n}\right)}{\sum_{i=1}^{n}w^{\ast}_{i-1}}\right)\\
& = & \sqrt{h_{n} n \Delta_n }\left(\frac{\delta_n}{\frac{1}{n^2
h_{n}^{4}}\sum_{i=1}^{n}w^\ast_{i-1}}\right)+o_p(1).
\end{eqnarray*}
Using the Lemma \ref{l3}, we have $\frac{1}{n^2
h_{n}^{4}}\sum_{i=1}^{n}w^\ast_{i-1} \stackrel{p}{\rightarrow}
p^{2}(x)[K_{1}^{2} - (K_{1}^{1})^{2}].$ Hence,
$$\hat{U}_n(x)- U^\ast_n(x)=\sqrt{h_{n} n \Delta_n }O_p(\Delta_n) =
O_{p}(\sqrt{h_{n} n \Delta_n^{3} }) \stackrel{p}{\rightarrow} 0$$ by
Assumption \ref{a5}.
\end{proof}
\bigskip

{\bf Acknowledgments} This research work is supported by Ministry of
Education, Humanities and Social Sciences project (No. 17YJA790075),
the General Research Fund of Shanghai Normal University (No.
SK201720) and Funding Programs for Youth Teachers of Shanghai
Colleges and Universities (No. A-9103-17-041301).

\bibliographystyle{amsplain}

\end{document}